\documentclass{amsart}
\usepackage{amsmath}
\usepackage{amsfonts}
\usepackage{mathrsfs}
\usepackage{bbm}
\usepackage{bm}
\usepackage{dsfont}
\usepackage[normalem]{ulem}
\usepackage{enumerate}
\usepackage{xmpmulti}
\usepackage[applemac]{inputenc}
\usepackage[english]{babel}
\usepackage[T1]{fontenc}
\usepackage{amsthm, overpic}
\usepackage{verbatim}
\usepackage{tikz}

\usepackage%
{hyperref}

\newcommand{\bpf}[1][Proof]{{\noindent {\sc #1: }}}
\newcommand{\epf}{{{\hfill $\Box$ \smallskip}}}
\newcommand{\IPM}{invariant probability measure} 
\newcommand{\CCs}{connected components} 
\newcommand{\R}{\mathbb{R}}
\newcommand{\N}{\mathbb{N}}

\newcommand{\tbf}{\mathbf{t}}

\newcommand{\ibf}{\mathbf{i}}

\newcommand{\Pp}{\mathsf{P}}

\newcommand{\E}{\mathbf{E}}
\newcommand{\Ic}{\mathcal{I}}

\newcommand{\Bc}{\mathcal{B}}

\newcommand{\AC}{\mathcal{A}}
\newcommand{\M}{\mathcal{M}}

\newcommand{\eps}{\epsilon}

\newcommand{\ONE}{\mathds{1}}

\newcommand{\Qc}{\mathcal{Q}}
\newcommand{\id}{\mathbbm{1}}
\newcommand{\Leb}{\text{Leb}} 

\newtheorem{theorem}{Theorem}
\newtheorem{lemma}{Lemma}
\newtheorem{remark}{Remark}
\newtheorem{corollary}{Corollary}
\newtheorem{proposition}{Proposition}

\title[Singularities of invariant densities]{Singularities of invariant densities for random switching between two linear ODEs in 2D.}

\author{Yuri Bakhtin, Tobias Hurth, Sean D. Lawley, Jonathan C.
  Mattingly} \address{Courant Institute of Mathematical Sciences, New
  York University, 251 Mercer St, New York, NY 10012 USA}
\address{Universit\'e de Neuch\^atel, Institut de math\'ematiques, Rue Emile-Argand 11, CH-2000 Neuch\^atel}
\address{Department of Mathematics, University of Utah, Salt Lake
  City, UT 84112 USA} 
\address{Mathematics Department and Department of Statistical Science, Duke
  University, Durham, NC 27708 USA}

\begin{document}
\begin{abstract}

We consider a planar dynamical system generated by two stable linear
vector fields with distinct fixed points and random switching between
them. We characterize singularities of the invariant density in terms
of the switching rates and contraction rates. We prove boundedness
away from those singularities. We also discuss some motivating
biological examples.

\end{abstract}

\maketitle

\section{Introduction}      \label{sec:introduction}

This paper describes the formation of singularities and regularity properties
 in the stationary
densities for the dynamics created by random switching between two
linear ordinary differential equations (ODEs) in the 
two-dimensional plane. 
%
%
A full characterization of stationary density singularities for 
randomly switched ODEs in one dimension is provided in \cite{Mattingly}, yet
singularity formation  is poorly understood in higher dimensions,
even at the level of motivational examples.

Here, we study a deceptively simple two-dimensional example in
the hope that it will begin to illuminate a path forward.  We have
not sought generality; but rather, picked a simple switching system between two linear equations
 to explore
how geometry of contraction and random switching interact to produce
singularities in the longtime distribution of the system. 
Despite this apparent
simplicity, the structure of the stationary density can be quite
rich. 
Depending on the
relationships between the switching rates and the contraction
rates, the stationary density may be bounded, have isolated singularites,  or have one-dimensional curves of singularities.
  Though we have studied a particular system, our methods are fairly
general and hopefully can be extended to an interesting class of examples.

There has been a resurgence in the study of such switched ODE systems
in recent years under the names \emph{hybrid systems} \cite{Yin}, \emph{piecewise
deterministic Markov processes (PDMP)} \cite{Davis, Malrieu_2015}, and
\emph{random evolutions}~\cite{Hersh-story:MR1962927}. Some of this renewed
interest stems from applications in ecology and cellular biology
\cite{LawleyMattinglyReed2015, bressloff2017rev}.
On the more theoretical side, 
it
was shown in~\cite{Bakhtin, Benaim, Malrieu} that a combination of a
condition of H\"ormander type and an accessibility condition
guarantees that an invariant distribution, if it exists, is unique and
absolutely continuous with respect to the Lebesgue measure.

One could expect that, similarly to the well-known results for
hypoelliptic diffusions based on pseudo-differential calculus or
Malliavin calculus, the same H\"ormander condition would guarantee
$C^\infty$ smoothness of the invariant density if the driving vector
fields are smooth and a hypoellipticiy condition is met. 
As already alluded to, the picture is more involved and
invariant densities of switching systems often have singularities. In \cite{Mattingly}, emergence of 
singularities of  invariant densities for one-dimensional switching systems due to contraction near stable
critical points was studied,  and
a classification of singularities was given. It was also shown that away from critical points of the driving vector
fields, the invariant densities are $C^\infty$. 

In higher dimensions, the situation is even more complex generically. Some of the flows generated by the driving vector fields may exhibit  long-term contraction with or without convergence to a stable critical point, e.g.,
 there may be more sophisticated low-dimensional
attractors.
Density singularities created by some of the vector fields may be propagated in new directions by other vector fields. Additional complexity emerges due to the presence of manifolds of hypoellipticity points.

We started an exploration of higher dimensions in \cite{torus}, where we
considered a class of switching systems on the two-dimensional torus that is devoid of these obstacles (the contraction is subexponential and all points are elliptic). For this class, we showed that the invariant densities belong to $C^\infty$ and that there are no singularities.

For generic switching systems, characterizing singularities of the invariant densities and proving smoothness away from those singularities still seems to be a hard problem. 
In the present paper, for the first time we consider switching systems with a whole line of points of hypoellipticity and contractive flows associated to the driving vector fields.  In various regimes that we define in terms of the parameters of the model, i.e., contraction rates and switching rates, we describe points and lines of singularities of the invariant density and prove
boundedness of the density away from those singularities. 

\medskip

Let us describe the system more precisely now. We consider the PDMP given by Poissonian random switching between the linear vector fields
\begin{equation}    \label{eq:switching_system} 
u_i(x_1, x_2) = \begin{pmatrix}  
                           -\alpha & 0 \\
                           0 & -\beta 
                         \end{pmatrix}  \begin{pmatrix} 
                                               x_1 - i \\
                                               x_2 - i 
                                              \end{pmatrix}, \quad i = 0, 1,  
\end{equation} 
where $\alpha > \beta > 0$. 
Given a starting point $x \in \R^2$ and an initial vector field, say $u_0$, we follow the flow of $u_0$ for an exponential time. Then a switch occurs, meaning that the driving vector field $u_0$ is replaced with $u_1$. Starting from the point in $\R^2$ where the switch occurred, we flow along $u_1$ for another exponential time, then switch back to $u_0$, etc.  We assume that the times between consecutive switches are independent. Switches from $u_0$ to $u_1$ happen at a constant rate $\lambda_0 > 0$, and switches from $u_1$ to $u_0$ happen at a constant rate $\lambda_1 > 0$.  The resulting dynamics are strongly affected by the globally asymptotically stable equilibrium points $(0,0)$ and $(1,1)$ of the two vector fields:  A typical switching trajectory obtained from intermittent switching between $u_0$ and $u_1$ enters in finite time the region $\Gamma$ bounded by the trajectory of $u_0$ starting from $(1,1)$ and the trajectory of $u_1$ starting from $(0,0)$, and then remains in $\Gamma$ for all future times (see Figure~\ref{fig:gamma} below).  Since the setting is essentially compact, the semigroup of the PDMP admits an invariant probability measure.  As will be established rigorously in Proposition~\ref{prop:support}, the invariant probability measure is unique and has a density with respect to Lebesgue measure. The goal of this article is to investigate the marginals $\rho_0$ and $\rho_1$ of the density, corresponding to the driving vector fields $u_0$ and $u_1$.  In this introduction and throughout the paper, we use the term \emph{invariant densities} for the marginals of the density associated with an absolutely continuous invariant probability measure. 

The PDMP governed by $u_0$ and $u_1$ can be thought of as a two-dimensional version of one of the simplest possible switching systems on the real line:  If we switch between $v_0(x) = -a x$ and $v_1(x) = a (1-x)$ for $a > 0$, the resulting switching trajectory is alternately attracted by $0$ and $1$.  As in the more complex two-dimensional system in~\eqref{eq:switching_system}, this simple one-dimensional system gives rise to a unique and absolutely continuous invariant probability measure.  Unlike the invariant densities in the 2D system, however, the invariant densities in the 1D system can be computed explicitly by solving the corresponding Kolmogorov forward equations, see e.g.~\cite{Faggionato}. They are densities of beta distributions: 
\begin{align}    
\begin{split}\label{eq:explicit_density_1}
\rho_0(x) =& c_0 x^{\frac{\lambda_0}{a}-1} (1-x)^{\frac{\lambda_1}{a}},  \\
\rho_1(x) =& c_1 x^{\frac{\lambda_0}{a}} (1-x)^{\frac{\lambda_1}{a}-1},   
\end{split}
\end{align}
where $c_0, c_1$ are constants.  In particular, $\rho_0$ and $\rho_1$ are smooth in the interior of $[0,1]$, and develop singularities at the critical points $0$ and $1$ if the switching rates are small compared to the rate of contraction $a$.  While it is possible to write down the Kolmogorov forward equations for the invariant densities of~\eqref{eq:switching_system}, we cannot find explicit solutions to the equations. Besides, it is a priori not clear whether the invariant densities of~\eqref{eq:switching_system} are sufficiently regular to be classical solutions on some meaningful set, say in the interior of $\Gamma$.  Notice, however, that the marginals of the invariant densities with respect to the coordinates $x_1$ and $x_2$ are explicitly given by the formulae in~\eqref{eq:explicit_density_1} for $a=\alpha$ and $a=\beta$.   We conjecture that the invariant densities for~\eqref{eq:switching_system} are $C^{\infty}$ in the interior of the set $\Gamma$. At the boundary of $\Gamma$, singularities may form due to exponential contraction and thus accumulation of probabilistic mass near the critical points, and the subsequent propagation of mass along trajectories of $u_0$ and $u_1$.  

We give two results on singularities of the invariant densities for slow switching. The first one describes the singularities near the  attracting critical points of $u_0$ and $u_1$.  The basic mechanism leading to these singularities is mass accumulation due to the fact that, under small to moderate switching rates, there are long time intervals during which the system is exposed to contraction towards $(0,0)$ and $(1,1)$. The second result holds only for small switching rates, and describes how a singularity at the critical point of $u_1$ is spread along the trajectory of $u_0$ passing through this critical point.  These results on singularity formation are complemented by several boundedness results, a first step towards proving regularity of the invariant densities in the interior of $\Gamma$. For instance, we show that the invariant densities are bounded on any compact set contained in the interior of $\Gamma$, even if switches are rare. 

There are two main difficulties in dealing with the switching system in~\eqref{eq:switching_system}.  The main obstacle to showing smoothness of the invariant densities is arguably the exponential contraction in the vicinity of critical points. Another less obvious difficulty stems from the fact that the vector fields $u_0$ and $u_1$ are aligned with each other along the diagonal line segment connecting $(0,0)$ and $(1,1)$. This partial breakdown of transversality makes the smoothing effect of switches close to the diagonal less pronounced. On the other hand, switches close to the diagonal but far from the critical points at least do not spoil the densities, which seems to make them a technical nuisance rather than an essential obstacle to establishing smoothness.

The paper is organized as follows. 
In Section~\ref{sec:applications}, we discuss two systems emerging in applications that can be reduced to \eqref{eq:switching_system}. We state our results on singularities of the invariant density in Section~\ref{sec:switching_system}. In Section~\ref{sec:support}, we prove existence and uniqueness of the invariant distribution, as well as a basic description of its support. Furthermore, we exhibit the line of hypoellipticity points, which is an obstacle to establishing boundedness of the invariant density. In Section~\ref{sec:integral_equation}, we recall some basic integral equations satisfied by the invariant density. In Section ~\ref{sec:proof_unbounded}, we prove one of our main results (Theorem~\ref{thm:unbounded}), which describes the singularities of the invariant density.
In Section~\ref{sec:change_of_variables}, we perform a change of variables in the integral equations from Section~\ref{sec:integral_equation} that prepares the proof of our main boundedness result (Theorem~\ref{thm:bounded}). The latter
is given in Section~\ref{sec:proof_bounded}, where most of the technical work is carried out.


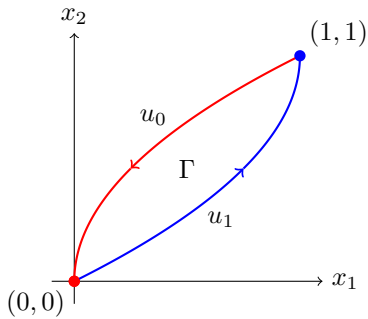
\begin{figure}\label{fig:gamma}
\begin{tikzpicture}[scale=3]
      \draw[->] (-.1,0) -- (1.1,0) node[right] {$x_{1}$};
      \draw[->] (0,-.1) -- (0,1.1) node[above] {$x_{2}$};
      \node[below left] at (0,0) {$(0,0)$};
      \node[above right] at (1,1) {$(1,1)$};
      \fill[color=red] (0,0) circle (.025);
  \draw[thick,scale=1,domain=0:.5,smooth,variable=\y,blue,->] plot ({2*\y-\y*\y},{\y});
    \draw[thick,scale=1,domain=.5:1,smooth,variable=\y,blue] plot ({2*\y-\y*\y},{\y});
     \draw[thick,scale=1,domain=1:.5,smooth,variable=\y,red,->]  plot ({\y*\y},{\y});
          \draw[thick,scale=1,domain=.5:0,smooth,variable=\y,red]  plot ({\y*\y},{\y});
         
          \node[above] at (0.35,0.65) {$u_{0}$};
          \node[below] at (0.65,0.35) {$u_{1}$};
          \node at (0.5,0.5) {$\Gamma$};
     
      \fill[color=red] (0,0) circle (.025);
      \fill[color=blue] (1,1) circle (.025);
    \end{tikzpicture}
    \caption{The support $\Gamma$ of the invariant densities is the region bounded by the forward $u_0$ trajectory starting at $(1,1)$ and the forward $u_1$ trajectory starting at $(0,0)$. 
}
\end{figure}

\section{Applications}
\label{sec:applications}

Generically, our results concern any two-dimensional randomly switching ODE of the form
\begin{align}\label{eq:gen}
\frac{d}{dt}\begin{pmatrix}
x_{1}\\ x_{2}
\end{pmatrix}
=A
\begin{pmatrix}
x_{1}\\ x_{2}
\end{pmatrix}
+\textbf{b}_{0} \mathbbm{1}_{I_t=0}+\textbf{b}_{1} \mathbbm{1}_{I_t=1},\quad \textbf{b}_{0},\textbf{b}_{1}\in\R^{2},
\end{align}
where $I_t \in\{0,1\}$ is a Markov jump process 
and $A\in\R^{2\times2}$ has two distinct, negative eigenvalues. In particular, \eqref{eq:gen} reduces to~\eqref{eq:switching_system} after the coordinate change 
$$
\begin{pmatrix}
y_1 \\
y_2 
\end{pmatrix} = G \begin{pmatrix} 
                                       x_1 \\
                                       x_2 
                                       \end{pmatrix} + \begin{pmatrix} 
                                                                - \frac{1}{\alpha} & 0 \\
                                                                 0 & - \frac{1}{\beta} 
                                                                 \end{pmatrix} G \textbf{b}_0 \mathbbm{1}_{I_t = 0} + \left(\begin{pmatrix} 
                                                                                                                                                         -\frac{1}{\alpha} & 0 \\
                                                                                                                                                         0 & -\frac{1}{\beta} 
                                                                                                                                                         \end{pmatrix} G \textbf{b}_1 + \begin{pmatrix} 
                                                                                                                                                                                                          1 \\
                                                                                                                                                                                                           1 
                                                                                                                                                                                                           \end{pmatrix} \right) \mathbbm{1}_{I_t = 1}, 
                                                                                                                                                                                                           $$
where $G$ is an invertible $(2 \times 2)$ matrix such that 
$$
A = G^{-1} \begin{pmatrix}
                  - \alpha & 0 \\
                  0 & -\beta 
                  \end{pmatrix} G. 
$$
In addition to being one of the simplest nontrivial two-dimensional PDMP examples in which to study invariant densities, models of the form~\eqref{eq:gen} arise naturally in diverse applications. We now give two such applications.

\subsection{Stochastic gene expression}

Much of the recent interest in PDMPs stems from their application to gene expression \cite{kepler2001,paulsson2005,bressloff2017rev}. Models in this context typically begin with a continuous-time Markov chain on a discrete state space that tracks gene products (an integer number of mRNA and/or protein molecules) as well as some discrete (often binary) environmental state, such as whether or not a gene is active or inactive. Assuming that the number of gene products is large, one often approximates the amount of gene product by a continuous variable that evolves by a deterministic ODE between stochastic switches in the environmental state. That is, the stochasticity stemming from the finite number of gene products is averaged out, while the stochastic environmental state is retained. 

To illustrate this concretely, we briefly describe the so-called ``standard model'' of gene expression \cite{paulsson2005}. Let $I_t \in\{0,1\}$ be the state of a gene, with $I_t=0$ ($I_t=1$) corresponding to an active (inactive) gene, and suppose $I_t$ leaves state $i\in\{0,1\}$ at rate $\lambda_i >0$. When the gene is active, it produces mRNA molecules at rate $\alpha>0$. Each mRNA molecule degrades at rate $\delta>0$ and produces a protein molecule at rate $\beta>0$. Protein molecules degrade at rate $\gamma>0$. Letting $X_t \in\{0,1,2,\dots\}$ and $Y_t \in\{0,1,2,\dots\}$ denote the respective mRNA and protein copy numbers, the Markov transitions are summarized by
\begin{align}\label{eq:mc}
\begin{split}
&I_t:\quad0\underset{\lambda_{1}}{\overset{\lambda_{0}}{\rightleftharpoons}}1;\qquad
X_t:\quad X\underset{\delta}{\overset{I_t \alpha}{\rightleftharpoons}}X+1,\qquad
Y_t:\quad Y\underset{\gamma}{\overset{X_t \beta}{\rightleftharpoons}}Y+1.
\end{split}
\end{align}
This three-component Markov chain $(X_t,Y_t,I_t)\in\{0,1,2,\dots\}^{2}\times\{0,1\}$ and various simplifications have been very well studied using a variety of mathematical techniques \cite{paulsson2005,bressloff2017rev}. Indeed, depending on the parameter regime, this Markov chain has been reduced to an ODE, a PDMP, a stochastic differential equation (SDE) driven by white noise, an SDE driven by L\'{e}vy noise, and a  L\'{e}vy-type process \cite{jia2017}.

For our purposes, suppose that the characteristic number of mRNA and protein molecules is large,
\begin{align*}
X^{*}:=\frac{\alpha}{\delta}\frac{\lambda_{0}}{\lambda_{0}+\lambda_{1}}\gg1,\quad
Y^{*}:=\frac{\beta}{\gamma}X^{*}\gg1.
\end{align*}
In this parameter regime, one can approximate the rescaled mRNA and protein concentrations, $x(t):=X_t/X^{*}$ and $y(t):=Y_t/Y^{*}$, by the two-dimensional PDMP \cite{smiley2010},
\begin{align}
\begin{split}\label{eq:xy}
\tfrac{d}{dt}x(t)
&=\frac{\alpha}{X^{*}} I_t -\delta x(t),\\
\tfrac{d}{dt}y(t)
&=\gamma (x(t)-y(t)),
\end{split}
\end{align}
in which the only source of stochasticity remaining is $I_t$. Of course, \eqref{eq:xy} is of the form~\eqref{eq:gen}.

\subsection{PDEs with randomly switching boundary conditions}

While most of the interest in PDMPs has focused on switching ODEs, a number of biological applications have recently prompted the study of PDEs with randomly switching boundary conditions (for example, see \cite{LawleyMattinglyReed2015,lawley2018,bressloff2017dyn,lawley2016}). Perhaps the simplest such example is the one-dimensional diffusion equation,
\begin{align*}
\tfrac{\partial}{\partial t}c(x,t)=\tfrac{\partial^{2}}{\partial x^{2}}c(x,t),\quad x\in(0,1),
\end{align*}
with an absorbing boundary condition at $x=0$ and a randomly switching boundary condition at $x=1$,
\begin{align*}
c(0,t)=0,\quad
c(1,t)= I_t,
\end{align*}
where $I_t\in\{0,1\}$ is a continuous-time Markov jump process. Writing the solution in terms of the $L^{2}[0,1]$-orthonormal basis, $\{\sqrt{2}\sin(n\pi x)\}_{n=1}^{\infty}$,
\begin{align*}
c(x,t)
=\sum_{n=1}^{\infty}c_{n}(t)\sqrt{2}\sin(n\pi x),
\end{align*}
it follows that any pair of coefficients, say $c_{k}(t)$ and $c_{m}(t)$, satisfy the two-dimensional switching ODEs,
\begin{align}\label{eq:modes}
\begin{split}
\tfrac{d}{dt}c_{k}(t)
&=-\beta_{k}(c_{k}(t)- I_t b_{k}),\\
\tfrac{d}{dt}c_{m}(t)
&=-\beta_{m}(c_{m}(t)- I_t b_{m}),
\end{split}
\end{align}
where $\beta_{n}=n^{2}\pi^{2}$ and $b_{n}=(-1)^{n+1}\sqrt{2}/(n\pi)$. Of course, \eqref{eq:modes} is of the form~\eqref{eq:gen}.

\bigskip

\section{Problem setting and main results}    \label{sec:switching_system}

We consider random switching between the linear vector fields $u_0$ and $u_1$ on $\R^2$, given by 
\begin{equation*}
u_i(x) =  u_i(x_1,x_2) = \begin{pmatrix}
           -\alpha & 0 \\
           0 & -\beta 
          \end{pmatrix} \begin{pmatrix}
                                   x_1 - i \\
                                   x_2 - i
                               \end{pmatrix}, \quad i = 0,1,   
\end{equation*}
where $\alpha > \beta > 0$. 
The vector fields $u_0$ and $u_1$ have an attracting critical point at $(0,0)$ and $(1,1)$, respectively.  For any
$(x_1, x_2) \in \R^2$, the initial-value problem
\begin{equation*}
 \begin{pmatrix}
  \dot{x_1}(t) \\
  \dot{x_2}(t)
 \end{pmatrix} = u_i(x_1(t), x_2(t)), \quad \begin{pmatrix}
                                           x_1(0) \\
                                           x_2(0)
                                        \end{pmatrix} = \begin{pmatrix}
                                                            x_1 \\
                                                            x_2
                                                         \end{pmatrix}
\end{equation*}
has the unique solution 
\begin{equation}    \label{eq:phi_v}
 \Phi_i^t(x_1, x_2) = \begin{pmatrix}
                     i + (x_1 -i) e^{-\alpha t} \\
                     i + (x_2 - i) e^{-\beta t}
                 \end{pmatrix}, \quad t \in \R. 
\end{equation}
It is easy to see that 
\begin{equation}       \label{eq:u_and_v}
 \Phi_1^t(x_1, x_2) = \begin{pmatrix}
                        1 \\
                        1
                       \end{pmatrix} - \Phi_0^t(1 - x_1, 1 - x_2), \quad
x \in \R^2, \ t \in \R. 
\end{equation}
For notational convenience, we also define the inverse flows 
$$ 
\Psi_i^t(x) = (\Phi_i^t)^{-1}(x) = \Phi_i^{-t}(x), \quad i \in \{0,1\}, \ t \in \R, \ x \in \R^2. 
$$ 
As we will be switching intermittently between $u_0$ and $u_1$, it is also convenient to define the cumulative flows 
$$ 
\Phi_i^{(t_1, \ldots, t_n)} = \begin{cases}
                                                 \Phi_i^{t_n} \circ \Phi_{1-i}^{t_{n-1}} \circ  \Phi_i^{t_{n-2}} \circ \ldots \circ \Phi_{1-i}^{t_1}, & \quad n \equiv 0 \mod 2, \\
                                                \Phi_i^{t_n} \circ \Phi_{1-i}^{t_{n-1}} \circ  \Phi_i^{t_{n-2}} \circ \ldots \circ \Phi_i^{t_1}, & \quad n \equiv 1 \mod 2
\end{cases} 
$$ 
and 
$$ 
\Psi_i^{(t_1, \ldots, t_n)} = \left( \Phi_i^{(t_1, \ldots, t_n)} \right)^{-1}. 
$$ 
For $i \in \{0,1\}$, we call the set $\{ \Phi_i^t(x): \ t > 0\}$ the forward $u_i$ trajectory
starting at $x$ and we call $\{ \Phi_i^t(x): \ t < 0\}$
the backward $u_i$ trajectory starting at $x$. The set $\{ \Phi_i^t(x): \ t
\in \R\}$ is simply called the $u_i$ trajectory through $x$. 

Let $I = (I_t)_{t \geq 0}$ be a continuous-time Markov chain on $\{0,1\}$ with jump rate $\lambda_0$ from $0$ to $1$ and $\lambda_1$ from $1$ to $0$. Then, we define a stochastic process $X = (X_t)_{t \geq 0}$ on $\R^2$ via 
\begin{equation}
  \label{eq:Xt}
  \tfrac{d}{dt} X_t = u_{I_t}(X_t). 
\end{equation}
The two-component process $(X,I)$ is a Markov process on $\R^2 \times
\{0,1\}$, whose Markov semigroup we denote by $(\Pp^t)_{t \geq 0}$. We
call a probability measure $\mu$ on $\R^2 \times \{0,1\}$ an
\emph{invariant probability measure} 
of $(\Pp^t)_{t \geq 0}$ if $\mu = \mu \Pp^t$ for all $t \geq 0$. 

The forward $u_0$ trajectory starting at $(1,1)$ and the forward
$u_1$ trajectory starting at $(0,0)$ together with the critical points $(0,0)$ and $(1,1)$ mark the boundary of
the set 
\begin{equation*}
 \Gamma = \left\{ (x_1,x_2) \in \R^2:  0 \leq x_2 \leq 1, \ x_2^{\frac{\alpha}{\beta}} \leq x_1 \leq 1 - (1-x_2)^{\frac{\alpha}{\beta}}\right\}. 
\end{equation*}
We denote the interior of $\Gamma$ by $\Gamma^{\circ}$. Notice that $\Gamma$ and $\Gamma^{\circ}$ are symmetric about the point $(\tfrac{1}{2},
\tfrac{1}{2})$, i.e. $(x_1,x_2) \in \Gamma$ ($\in \Gamma^{\circ}$) if and only if $(1 - x_1, 1 -
x_2) \in \Gamma$ ($\in \Gamma^{\circ}$). 

\begin{proposition} \label{prop:support} The Markov semigroup
  $(\Pp^t)_{t \geq 0}$ admits a unique {\IPM} $\mu$. It is absolutely
  continuous with respect to the product of Lebesgue measure on $\R^2$
  and counting measure on $\{0,1\}$. Moreover, the marginals
  $\mu_i(\cdot) = \mu(\cdot \times \{i\})$, $i \in \{0,1\}$, have
  support $\Gamma$.
\end{proposition}

Recall that the support of $\mu_i$ is the collection of all points $x \in \R^2$ such that $\mu_i(U) > 0$ for every neighborhood $U$ of $x$.  We prove Proposition~\ref{prop:support} in Section~\ref{sec:support}. 
Since the marginal $\mu_i$, $i \in \{0,1\}$, is absolutely continuous with respect to Lebesgue measure, it has a density $\rho_i \in L^1(\R^2)$, which we call an \emph{invariant density}. Below, we state our results on boundedness as well as the occurrence of singularities for the invariant density $\rho_0$.  Exploiting the symmetries of the
switching system, one can easily formulate corresponding results for
$\rho_1$.  Whether and where singularities of $\rho_0$ occur depends critically on the switching rates $\lambda_0$ and $\lambda_1$.  In some sense this is not surprising because small switching rates translate into few switches and thus an accumulation of probabilistic mass at the critical points $(0,0)$ and $(1,1)$. Interestingly, if both $\lambda_0$ and $\lambda_1$ are very small, the singularity created at the critical point of one of the vector fields is propagated along the forward trajectory of the other vector field that starts at the critical point.  As $L^1$ functions, $\rho_0$ and $\rho_1$ are only defined up to a set of Lebesgue measure zero, so when we state, e.g., that $\rho_0$ is bounded on a set $S$, we mean that there is a representative of $\rho_0$ that is bounded on $S$.  Proposition~\ref{prop:support} implies that $\rho_0$ and $\rho_1$ vanish outside of $\Gamma$, which is why we can restrict ourselves to $\Gamma^{\circ}$ instead of considering all of $\R^2$.  

For $i \in \{0,1\}$, let $\partial \Gamma_i$ denote the forward $u_i$ trajectory starting at $(1-i, 1-i)$, and set  
$$
\Gamma_i = \Gamma \setminus \{(i,i)\}. 
$$ 
Observe that $\partial \Gamma_i$, $i \in \{0,1\}$, are the curves that make up the right and left part of the boundary of $\Gamma$, minus the critical points $(0,0)$ and $(1,1)$. The following theorem describes for which switching rates and in which regions singularities occur. 

\begin{theorem}     \label{thm:unbounded}
The following statements hold. 
\begin{enumerate}
\item For $\lambda_0 < \alpha + \beta$, the invariant density $\rho_0$ is unbounded in every neighborhood of $(0,0)$. 
\item For $\lambda_1 < \beta$ and $x \in \partial \Gamma_0$, $\rho_0$ is unbounded in every neighborhood of $x$.  Since being unbounded in every neighborhood of a point is a closed condition, $\rho_0$ is also unbounded in every neighborhood of $(0,0)$ and $(1,1)$. 
\end{enumerate}
\end{theorem}

\begin{theorem}       \label{thm:bounded}
The following statements hold. 
\begin{enumerate}
\item For $\lambda_0 > \alpha + \beta$ and $\lambda_1 > \beta$, the invariant density $\rho_0$ is bounded on $\Gamma^{\circ}$.
\item Let $\lambda_0, \lambda_1 > \beta$ and let $K \subset \Gamma_0$ be compact.  Then, $\rho_0$ is bounded on $K$. 
\item Let $K \subset \Gamma$ be a compact set such that $K \cap \partial \Gamma_0 = \emptyset$. Then, $\rho_0$ is bounded on $K$ for any switching rates $\lambda_0, \lambda_1 > 0$. 
\end{enumerate}
\end{theorem}

\begin{remark} \rm 
Theorems~\ref{thm:unbounded} and~\ref{thm:bounded} do not address whether $\rho_0$ stays bounded along $\partial \Gamma_0$ if $\lambda_0 < \beta$ and $\lambda_1 > \beta$. Based on simulations and heuristics, we conjecture that $\rho_0$ is bounded in this case, i.e. we conjecture the conclusion of Theorem~\ref{thm:bounded}, part (2), to hold for every $\lambda_0 > 0$ and $\lambda_1 > \beta$. The critical cases not covered by Theorems~\ref{thm:unbounded} and~\ref{thm:bounded} (e.g., is $\rho_0$ bounded on $\Gamma^{\circ}$ if $\lambda_0 = \alpha + \beta$ and $\lambda_1 = \beta$?) are also open.  
\end{remark}   

We prove Theorem~\ref{thm:unbounded} in Section~\ref{sec:proof_unbounded}. The proof of Theorem~\ref{thm:bounded} is given in Section~\ref{sec:proof_bounded}.   Theorems~\ref{thm:unbounded} and~\ref{thm:bounded} combined provide the following picture:  For fast switching away from $u_0$ ($\lambda_0 > \alpha + \beta$) and for intermediate switching away from $u_1$ ($\lambda_1 > \beta$), the invariant density $\rho_0$ is globally bounded. For intermediate switching away from $u_0$ ($\lambda_0 < \alpha + \beta$), $\rho_0$ has a singularity at $(0,0)$, the critical point of $u_0$, irrespective of $\lambda_1$. And in the regime of slow switching away from $u_1$ ($\lambda_1 < \beta$), $\rho_0$ has singularities along the entire left boundary curve of the support, including $(0,0)$. This happens regardless of how quickly on average we switch away from $u_0$. Away from the left boundary curve, $\rho_0$ is always bounded. The mechanism leading to the blow-up of $\rho_0$ along $\partial \Gamma_0$ can be roughly described as follows: Due to exponential contraction of the flow of $u_1$, probabilistic mass accumulates at the sink $(1,1)$. This mass is subsequently propagated under the flow of $u_0$ and thus gives rise to singularities along the forward $u_0$ trajectory starting at $(1,1)$.

\section{The support of the invariant measure}   \label{sec:support}

In this section, we prove Proposition~\ref{prop:support}.  Given two points $x, y \in \R^2$, we say that $x$ is \emph{reachable} from $y$ if there exist $i \in \{0,1\}$, $n \in \N$, and $(t_1, \ldots, t_n) \in \R^n_+$ such that 
$$ 
x = \Phi_i^{(t_1, \ldots, t_n)}(y). 
$$ 
For $y \in \R^2$, let $L(y)$ denote the set of points $x \in \R^2$ that are reachable from $y$. A nonempty set $S \subset \R^2$ is called \emph{positive invariant} if $L(y) \subset S$ for every $y \in S$. This is equivalent to saying that 
$$ 
\Phi_i^t(x) \in S, \quad i \in \{0,1\}, \ t \geq 0, \ x \in S. 
$$ 

\begin{lemma}   \label{lm:invariance}
The following statements hold. 
\begin{enumerate}
\item The sets $\Gamma$ and $\Gamma^{\circ}$ are positive invariant.
\item  We have $\Gamma^{\circ} \subset L(x)$ for every $x \in \R^2$.
\end{enumerate}  
\end{lemma}

\bpf To simplify notation, we define 
$$ 
\gamma := \frac{\alpha}{\beta} > 1. 
$$ 
Fix a point $x \in \Gamma$. If $x = (0,0)$, the forward $u_0$ trajectory starting at $x$ consists only of $x$ and is therefore contained in $\Gamma$. If $x \neq (0,0)$
and if $y$ is any point on the forward $u_0$ trajectory starting at $x$, we have 
\begin{equation*}
 y_2^{\gamma} \leq x_1 x_2^{-\gamma} y_2^{\gamma} \leq \left( 1 - (1-x_2)^{\gamma} \right) x_2^{-\gamma} y_2^{\gamma} = \left(\frac{y_2}{x_2} \right)^{\gamma} - \left( \frac{y_2}{x_2} - y_2 \right)^{\gamma} \leq 1 - (1-y_2)^{\gamma}. 
\end{equation*}
To obtain the last inequality, we used that 
$$ 
\frac{d}{dz} \left(z^{\gamma} - (z-y_2)^{\gamma} \right) = \gamma \left(z^{\gamma - 1} - (z-y_2)^{\gamma - 1} \right) > 0, \quad z > y_2, 
$$
and that $y_2 / x_2 < 1$. Since $x_1 x_2^{-\gamma} y_2^{\gamma} = y_1$, it follows that the forward
$u_0$ trajectory starting at $x$ is contained in $\Gamma$. 
As $(1 - x_1, 1 - x_2) \in \Gamma$, we also have 
\begin{equation}   \label{eq:invariance_2}
 \Phi_0^t(1 - x_1, 1 - x_2) \in \Gamma, \quad t > 0. 
\end{equation}
Equations~\eqref{eq:u_and_v} and~\eqref{eq:invariance_2} imply that the forward $u_1$ trajectory
starting at $x$ is contained in $\Gamma$ as well.  The proof for $\Gamma^{\circ}$ is analogous. 

\medskip

To prove the second statement, fix $x \in \R^2$ and $y \in \Gamma^{\circ}$. Since the critical points $(0,0)$ and $(1,1)$ are globally asymptotically stable, the set $L(x)$ contains a point $z$ such that $z_2 \in (0,1)$. If $z_1 \geq z_2$, set 
$$ 
a(t) = \Phi_0^t(z), \quad t > 0. 
$$ 
We have 
\begin{equation}      \label{eq:star_eq} 
\left( e^{-\beta t} z_2 \right)^{\gamma} < e^{-\alpha t} z_1. 
\end{equation}
Set 
$$ 
b(t) = e^{-\alpha t} z_1 + \left(1 - e^{-\beta t} z_2 \right)^{\gamma}. 
$$ 
Since $\lim_{t \to \infty} b(t) = 1$ and since 
$$ 
b'(t) = -\alpha z_1 e^{-\alpha t} + \alpha z_2 e^{-\beta t} \left(1 - e^{-\beta t} z_2 \right)^{\gamma-1} > 0
$$ 
for $t$ sufficiently large, we have $b(t) < 1$ for large $t$. 
This and~\eqref{eq:star_eq} imply that $a(t) \in \Gamma^{\circ}$ for large $t$.  In particular, there is $a \in L(x) \cap \Gamma^{\circ}$.  If $z_1 < z_2$, set 
$$ 
a(t) = \Phi_0^t(1-z_1, 1-z_2), \quad t > 0. 
$$ 
As $1 - z_1 > 1 - z_2$, we have as before $a(t) \in \Gamma^{\circ}$ for large $t$. With~\eqref{eq:u_and_v} and symmetry of $\Gamma^{\circ}$ around $(\tfrac{1}{2}, \tfrac{1}{2})$, this yields  
$$ 
\Phi_1^t(z) \in \Gamma^{\circ}
$$ 
for large $t$, so as in the case $z_1 \geq z_2$ there is $a \in L(x) \cap \Gamma^{\circ}$.  

As $L(a) \subset L(x)$, we may assume without loss of generality that $x \in \Gamma^{\circ}$. 
Now, it suffices to show that one of the following statements holds, as this means there are $s, t > 0$ and $i \in \{0,1\}$ such that $y = \Phi_i^{(s,t)}(x)$.
\begin{enumerate}[(a)]
 \item There is $\eta \in (0, \min\{x_2, y_2\}]$ such that 
\begin{equation*}
g(\eta) :=  1 - (1 - y_1) (1 - y_2)^{-\gamma} (1 - \eta)^{\gamma} - x_1 x_2^{-\gamma} \eta^{\gamma} = 0.
\end{equation*}
\item There is $\eta \in [\max\{x_2, y_2\}, 1)$ such that 
\begin{equation*}
h(\eta) := 1 - (1 - x_1) (1 - x_2)^{-\gamma} (1 - \eta)^{\gamma} - y_1 y_2^{-\gamma} \eta^{\gamma} = 0. 
\end{equation*}
\end{enumerate}
It is easy to see that $g(0)$, $g(1)$, $h(0)$, and $h(1)$ are negative. Since
$g$ and $h$ are continuous, it is then enough to show that $g(\min\{x_2,
y_2\}) \geq 0$ or $h(\max\{x_2, y_2\}) \geq 0$. First, assume that
$x_2 \leq y_2$. If 
$$
(1- y_1) (1 - y_2)^{-\gamma} \leq
(1- x_1) (1 - x_2)^{-\gamma},
$$
we have 
\begin{equation*}
 g(x_2) = 1 - (1 - y_1) (1 - y_2)^{-\gamma} (1 - x_2)^{\gamma} -
x_1 \geq 0.
\end{equation*}
And if 
$$
(1 - x_1) (1 - x_2)^{-\gamma} \leq (1 - y_1) (1 -
y_2)^{-\gamma},
$$
we have 
\begin{equation*}
 h(y_2) = 1 - (1 - x_1) (1 - x_2)^{-\gamma} (1 - y_2)^{\gamma} -
y_1 \geq 0. 
\end{equation*}
It remains to consider the case $x_2 > y_2$. If $x_1 x_2^{-\gamma} \leq
y_1 y_2^{-\gamma}$, we have $g(y_2) \geq 0$. And if
$y_1 y_2^{-\gamma} \leq x_1 x_2^{-\gamma}$, we have $h(x_2)
\geq 0$. 
\epf

\bigskip

Since the semigroup $(\Pp^t)_{t \geq 0}$ is Feller and since $\Gamma$ is compact and positive invariant, $(\Pp^t)_{t \geq 0}$ admits an \IPM. In fact, as we will show below, the {\IPM} is unique and absolutely continuous. For $x \in \R^2$, let 
$$ 
U(x) = (u_1(x), u_0(x))
$$ 
be the $(2 \times 2)$ matrix whose first column is $u_1(x)$ and whose second column is $u_0(x)$. As stated in the following easily verified lemma,  $u_0$ and $u_1$ are transversal at every point except for points on the line $x_1 = x_2$.

\begin{lemma}     \label{lm:uniqueness_ac}
Let $x \in \R^2$. Then, $\det U(x) = \alpha \beta (x_1 - x_2)$. In particular, $\det U(x) = 0$ if and only if $x_1 = x_2$.
\end{lemma}

In light of Lemmas~\ref{lm:invariance} and~\ref{lm:uniqueness_ac},
there are points (namely every point in $\Gamma^{\circ}$ not located
on $x_1 = x_2$) that are reachable from every starting point in $\R^2$
and where $u_0$ and $u_1$ are transversal.  By Theorem 1
in~\cite{Bakhtin} or by Theorem 4.5 in~\cite{Benaim}, the {\IPM} of
$(\Pp^t)_{t \geq 0}$ is unique and absolutely continuous with respect
to the product of Lebesgue measure on $\R^2$ and counting measure on
$\{0,1\}$.

Alternatively, one can follow the reasoning in~\cite{LawleyMattinglyReed2015}, which leverages the  contractive
nature of the system. This is particularly simple in this case as the
flows are deterministically uniformly contracting. Fixing two initial
condition $x,y \in \R^2$ and $i \in \{0,1\}$, we set $I_0=i$ and let
$X_t(x)$ and $X_t(y)$ be the solution to \eqref{eq:Xt} with the same
$I_t$ process (and hence the same jump times) but starting initially
from $x$ and $y$ respectively. 
If we define $r_t := X_t(x) - X_t(y)$, observe that
\begin{equation}
  \label{eq:contraction}
  \|r_t\| \leq \|x-y\| \exp( -(\alpha\wedge \beta) t),  
\end{equation}
where $\|\;\cdot\;\|$ denotes the Euclidean norm on $\R^2$. 
Now let $f\colon \R^2 \times \{0,1\}
\rightarrow \R$ be an arbitrary test function which is 1-Lipschitz continuous, i.e. $\lvert f(x,i) - f(y,j) \rvert \leq \| x - y \| + \id_{i \neq j}$ for every $x, y \in \R^2$ and $i, j \in \{0,1\}$. 
Since $\sup_{t \geq 0} \| X_t(x) \| < \infty$ for every $x \in \R^2$, the supremum of $(\Pp^tf)(x,i) - (\Pp^tf)(y,i) = \E f(X_t(x), I_t) -
\E f(X_t(y), I_t)$ over all such test functions is equal to the 
1-Wasserstein distance between $\Pp^t(x,i;\;\cdot\;)$ and
$\Pp^t(y,i;\;\cdot\;)$. Denoting this distance by $\| \Pp^t(x,i;\;\cdot\;) - \Pp^t(y,i;\;\cdot\;)\|_{W_1}$
and recalling that $\|\delta_{x,i} - \delta_{y,i}\|_{W_1} = \|x-y\|$ produces
\begin{align*}
 \|\Pp^t(x,i;\;\cdot\;)- \Pp^t(y,i;\;\cdot\;)\|_{W_1} \leq  \|\delta_{x,i} - \delta_{y,i}\|_{W_1}  e^{-(\alpha\wedge \beta) t}.
\end{align*}
A simple coupling argument using the definition of the 1-Wasserstein
distance as the infimum over all couplings and the resulting convexity of the 1-Wasserstein distance 
produces
\begin{align}   \label{eq:Wasserstein_estimate} 
 \|\mu \Pp^t- \nu \Pp^t\|_{W_1} \leq  \|\mu - \nu \|_{W_1}  e^{-(\alpha\wedge \beta) t}
\end{align}
for arbitrary initial probability measures $\mu$ and $\nu$ of bounded support. Observe in addition that there is a bounded positive invariant subset $B$ of $\R^2$ (e.g., the set $[-1,2]^2$) with the following property: For every $R > 0$ there is $C > 0$ such that every switching trajectory starting from a point of distance less than $R$ from the origin enters the set $B$ in a time less than $C$. As a result, every invariant probability measure has bounded support, so the estimate in~\eqref{eq:Wasserstein_estimate} proves in particular uniqueness of the invariant probability
measure. Additionally, this
proves a spectral gap for the Markov kernel $\Pp^t$ in the
1-Wasserstein distance which is independent of the switching
rates. The estimate in \eqref{eq:contraction} can also be used to show
that the system has a random attractor which consists of a  single
point and that all of the Lyapunov exponents are
negative. See~\cite{LawleyMattinglyReed2015} for more discussions in
this direction.

To finish the proof of Proposition~\ref{prop:support}, it remains to
show the statement about the support of the marginals $\mu_0$ and
$\mu_1$.  Since $\Gamma$ is positive invariant and compact, the
support of $\mu_0$ and $\mu_1$ is contained in $\Gamma$. And since
$\Gamma^{\circ}$ is an open subset of the set of points $y \in \Gamma$
that are reachable from every starting point in $\Gamma$,
$\Gamma^{\circ}$ is a subset of the support of $\mu_0$ and $\mu_1$,
see e.g.~\cite[Lemma 6]{Bakhtin}. As the support of $\mu_0$ and
$\mu_1$ is a closed set, it is necessarily equal to $\Gamma$.

\section{Integral equations for invariant densities and cdf's}        \label{sec:integral_equation}

Recall that for $i \in \{0,1\}$, $\Psi_i$ denotes the inverse flow associated with the vector field $u_i$.  Lemma 2 in~\cite{Mattingly} implies that 
\begin{equation}     \label{eq:pdf}
\rho_i(x) = \int_{\R_+} \lambda_{1-i} e^{-\lambda_i t} \det \nabla_x \Psi_i^t(x) \rho_{1-i}(\Psi_i^t(x)) \ dt, \quad i \in \{0,1\}.
\end{equation}
Written in terms of the cumulative distribution functions (CDF's)
\begin{equation}    \label{eq:cdf} 
G_i(x_1, x_2) = \int_{-\infty}^{x_1} \int_{-\infty}^{x_2} \rho_i(y_1, y_2) \ dy_1 \ dy_2, \quad (x_1, x_2) \in \R^2, 
\end{equation} 
the integral equations in~\eqref{eq:pdf} become 
\begin{equation}    \label{eq:G_0_as_integral_operator_of_G_1}
G_i(x) = \int_{\R_+} \lambda_{1-i} e^{-\lambda_i t} G_{1-i}(\Psi_i^t(x)) \ dt, \quad i \in \{0,1\}. 
\end{equation}
This is because for $i \in \{0,1\}$ and for any fixed $x \in \R^2$, $t \mapsto \Phi_i^t(x)$ is monotone in both components. Note that the integral on the right side of~\eqref{eq:G_0_as_integral_operator_of_G_1} can be rewritten as 
$$ 
\frac{\lambda_{1-i}}{\lambda_i} \E G_{1-i}(\Psi_i^T(x)), 
$$ 
where $T$ is an exponential random variable with intensity $\lambda_i$. 

Next, we generalize the integral equations in~\eqref{eq:pdf} by considering the evolution of $(X,I)$ leading to the current state not just since the latest switch but over the latest $n$ switches, $n \in \N$.  For $i \in \{0,1\}$, $n \in \N$, $\tbf \in \R^n_+$, and $x \in \Gamma^{\circ}$, we define the Jacobian 
$$ 
J_i^{\tbf}(x) = \det \nabla_x \Psi^{\tbf}_i(x). 
$$ 
For $i \in \{0,1\}$, $n \in \N$, and $x \in \Gamma^{\circ}$, let 
\begin{equation}
  \label{eq:TDef}
  T_i^n(x) = \left\{\tbf \in \R^n_+: \Psi_i^{\tbf}(x) \in \Gamma^{\circ} \right\}. 
\end{equation}
For $n \in \N$ and real-valued integrable functions $h$ on $\Gamma^{\circ}$, we define the transfer operator 
\begin{equation}  \label{eq:QcDef} 
\Qc_n h(x) = \int_{T^n_0(x)} \lambda_0(n) e^{-\langle \lambda_0^{(n)}, \tbf \rangle} J^{\tbf}_0(x) h(\Psi_0^{\tbf}x) \ d \tbf, \quad x \in \Gamma^{\circ},   
\end{equation}  
where 
\begin{equation}     \label{eq:lambda_def_1} 
\lambda_i(n) = \begin{cases}
                       \lambda_{1-i} \lambda_i \ldots \lambda_i, & \quad n \equiv 0 \mod 2, \\
                       \lambda_{1-i} \lambda_i \ldots \lambda_{1-i}, &  \quad n \equiv 1 \mod 2
\end{cases}
\end{equation}
is an alternating product of $\lambda_{1-i}$ and $\lambda_i$ with exactly $n$ factors, and where 
\begin{equation}      \label{eq:lambda_def_2} 
\lambda_i^{(n)} = \begin{cases}
                               (\lambda_{1-i}, \lambda_i, \ldots, \lambda_i)^{\top}, & \quad n \equiv 0 \mod 2, \\
                               (\lambda_i, \lambda_{1-i}, \ldots, \lambda_i)^{\top}, & \quad n \equiv 1 \mod 2
\end{cases}
\end{equation}
is a vector of length $n$ whose components alternate between $\lambda_i$ and $\lambda_{1-i}$. 
Then, 
\begin{equation}      \label{eq:Q_fixed_point}
 \rho_0 = \begin{cases}
                   \Qc_n \rho_0, & \quad n \equiv 0 \mod 2, \\
                   \Qc_n \rho_1, & \quad n \equiv 1 \mod 2, 
\end{cases}
\end{equation}
which can be deduced by iteratively plugging instances of~\eqref{eq:pdf} into one another and using the fact that the pushforward of a function under the cumulative flow $\Phi_i^{\tbf}$ is the composition of pushforwards under the individual flows $\Phi_i^{t_n}, \Phi_{1-i}^{t_{n-1}}, \ldots$

\begin{remark}   \rm 
The formula in~\eqref{eq:Q_fixed_point} can be generalized to switching systems with state space $U \times S$, where $U$ is an open subset of $\R^n$ and $S$ is a finite index set corresponding to a collection of smooth vector fields $u$ on $\R^n$ that leave $U$ positive invariant and are integrable, i.e. for any $x_0 \in U$ the initial-value problem $\dot{x} = u(x), x(0) = x_0$ has a unique solution, and this solution is defined for all $t \in \R$. Suppose the corresponding Markov semigroup admits an absolutely continuous invariant measure $\mu$ with invariant densities $(\rho_i)_{i \in S}$. For $i, j \in S$, let $\lambda_i$ be the rate of switching away from vector field $u_i$ and let $\lambda_{j,i}$ be the rate of switching from $u_j$ to $u_i$. For $n \in \N$, $\ibf = (i_1, \ldots, i_n) \in S^n$ and for real-valued integrable functions $h$ on $U$, define 
$$
\Qc_{\ibf} h(x) = \int_{T_{\ibf}(x)} \prod_{j=2}^n \lambda_{i_{j-1}, i_j} e^{-\sum_{j=1}^n \lambda_{i_j} t_j} \Phi_{\ibf}^{\tbf} \# h(x) \ d \tbf, \quad x \in U,
$$ 
where $\Phi_{\ibf}^{\tbf} = \Phi_{i_n}^{t_n} \circ \ldots \circ \Phi_{i_1}^{t_1}$, $T_{\ibf}(x) = \{\tbf \in \R^n_+: (\Phi_{\ibf}^{\tbf})^{-1}(x) \in U\}$ and where $\Phi_{\ibf}^{\tbf} \# h$ denotes the pushforward of $h$ under $\Phi_{\ibf}^{\tbf}$. Then, we have for $n \in \N$ and $i_n \in S$ that 
\begin{equation*}
\rho_{i_n} = \sum_{i_{n-1} \neq i_n} \ldots
\sum_{i_0 \neq i_1} \lambda_{i_0, i_1} \Qc_{(i_1, \ldots, i_n)} \rho_{i_0}. 
\end{equation*}
\end{remark}

\section{Singularities of the invariant densities}      \label{sec:proof_unbounded}

In this section, we prove Theorem~\ref{thm:unbounded} on singularities of the invariant density $\rho_0$ for slow switching.  We first show that if $\lambda_0 < \alpha + \beta$, then $\rho_0$ is unbounded in any neighborhood of $(0,0)$. 

\bigskip
\bpf[Proof of Theorem~\ref{thm:unbounded}, part (1)]  Recall that the cumulative distribution function $G_0$ was defined in~\eqref{eq:cdf}. Assuming that the invariant density $\rho_0$ is bounded by a constant $C$ in some neighborhood of $(0,0)$, we conclude that
\[G_0(\eps^{\alpha},\eps^{\beta})<C\eps^{\alpha+\beta}\] for sufficiently small $\eps$.
On the other hand, \eqref{eq:G_0_as_integral_operator_of_G_1} implies that, for every $\eps>0$,
\begin{align}    \label{eq:mass-near-origin}
 G_0(\eps^{\alpha},\eps^{\beta})&= \int_0^{\infty} \lambda_1 e^{-\lambda_0 t} G_1(e^{\alpha t}\eps^{\alpha}, e^{\beta t}\eps^{\beta}) \ dt  \\
                 &\ge  \int_{\ln(1/\eps)}^{\infty} \lambda_1 e^{-\lambda_0 t}G_1(e^{\alpha t}\eps^{\alpha}, e^{\beta t} \eps^{\beta}) \ dt.  \notag
\end{align}
For $t \geq \ln(\tfrac{1}{\epsilon})$, we have $e^{\beta t} \epsilon^{\beta} \geq 1$ and thus, as the support $\Gamma$ of $\mu_1$ is contained in $[0,1]^2$, we have $G_1(e^{\alpha t} \epsilon^{\alpha}, e^{\beta t} \epsilon^{\beta}) = G_1(+\infty,+\infty)$. Hence, the integral in the second line of~\eqref{eq:mass-near-origin} equals 
$$
G_1(+\infty, +\infty) \int_{\ln(1/\eps)}^{\infty} \lambda_1 e^{-\lambda_0 t}\ dt = G_1(+\infty, +\infty) \frac{\lambda_1}{\lambda_0} \eps^{\lambda_0}.
$$
Since $\eps$ is arbitrary, we conclude that $\lambda_0 \ge \alpha + \beta$. 
\epf

The idea behind the above proof is very general and can be
used with minor modifications to study existence and character of singularities in various other situations including 
high-dimensional ones. However, there is another interesting proof specific to the concrete vector fields $u_0,u_1$ we consider.

\bigskip 
  
\bpf[Another proof of Theorem~\ref{thm:unbounded}, part (1)]
Let $\epsilon\in(0,1)$. Since the point $(\epsilon^{\alpha},\epsilon^{\beta}) \in\R^2$ is on $\partial \Gamma_0$, the left boundary curve of $\Gamma$, and since $\rho_{0}$ is identically zero outside of $\Gamma$, we have 
\begin{align*}  
G_0(\eps^{\alpha},\eps^{\beta})=\int_{-\infty}^{\epsilon^{\alpha}}\int_{-\infty}^{\infty}\rho_0(y_1,y_2)\,dy_1\,dy_2=G_0(\eps^{\alpha},+\infty).
\end{align*}
Since our system can be viewed
as a product of non-iteracting components, the marginal distribution 
$$ 
E \mapsto \int_E \int_{-\infty}^{\infty} \rho_0(y_1, y_2) \ d y_1 \ d y_2
$$ 
coincides with the stationary distribution of the one-dimensional system given by switching at rates $\lambda_0$ and $\lambda_1$ between the one-dimensional vector fields 
$$ 
v_0(x) = -\alpha x, \quad v_1(x) = - \alpha (x-1). 
$$ 
By Proposition 3.12 in~\cite{Faggionato}, this distribution is a beta distribution with parameters $(\lambda_0 / \alpha,\lambda_1 / \alpha +1)$, so 
assuming that the invariant density $\rho_0$ is bounded by a constant $C$ in a neigborhood of $(0,0)$, we obtain for small $\eps$ 
\begin{align}\label{sup}
c B(\epsilon^{\alpha};\lambda_0/ \alpha,\lambda_1 / \alpha +1) = G_0(\eps^{\alpha},\eps^{\beta})\le C\epsilon^{\alpha + \beta},
\end{align}
where $c$ is a normalizing constant and where $B(x;a,b)$ is the incomplete beta function
\begin{align*}
B(x;a,b)=\int_0^x t^{a-1}\,(1-t)^{b-1}\,dt.
\end{align*}
Now if $\lambda_0 < \alpha + \beta$, it is straightforward to check that
\begin{align*}
\epsilon^{-(\alpha + \beta)} B(\epsilon^{\alpha};\lambda_0/ \alpha,\lambda_1/ \alpha +1) \to\infty,\quad\text{as }\epsilon\to0.
\end{align*}
Thus, dividing (\ref{sup}) by $\epsilon^{\alpha + \beta}$ and taking $\epsilon\to0$ completes the proof.
\epf

\bigskip

Now we prove the second part of Theorem~\ref{thm:unbounded}, which asserts that $\rho_0$ blows up along the entire left boundary curve of $\Gamma$ in case $\lambda_1 < \beta$.  

\bigskip

\bpf[Proof of Theorem~\ref{thm:unbounded}, part (2)]   Let us introduce the functions 
\begin{align*}
\phi^{(1)}(t)= e^{-\alpha t},\quad
\phi^{(2)}(t)=e^{- \beta t},\quad \phi_{\eps}^{(2)}(t)=(1-\eps^{\beta}) e^{- \beta t}, \quad t\in\R,\ \eps>0,  
\end{align*}
so that $\{(\phi^{(1)}(t),\phi^{(2)}(t))\}_{t > 0}$ is the forward $u_0$ trajectory starting at $(1,1)$ --- the left boundary curve of $\Gamma$.
Also, $\Phi_0^t(1,1-\eps^{\beta})=(\phi^{(1)}(t),\phi_{\eps}^{(2)}(t))$. 

For $z>0$, we define $t(z)= (\phi^{(1)})^{-1}(z)=-\ln(z)/ \alpha$ and the open interval $I_{\eps}(z)=(\phi_{\eps}^{(2)}(t(z)),\phi^{(2)}(t(z)))\subset\R$.
Also, for any set $I\subset \R$, we introduce
\[
R_{\eps}(I)=\left\{(x_1,x_2): x_1\in \phi^{(1)}(I),\  x_2\in I_{\epsilon}(x_1) \right\}. \]
Note that for all $t\in\R$ and all $t_1,t_2\in\R$ satisfying $t_1<t_2$, $\Phi_0^t:R_{\eps}(t_1,t_2)\to R_{\eps}(t_1+t,t_2+t)$ is a diffeomorphism. 
Also, for any $I,J\subset\R$,
\begin{equation}   \label{eq:R_identity}
R_{\eps}(I)\cap R_{\eps}(J)= R_{\eps}(I\cap J). 
\end{equation}

Suppose that there is $t > 0$ such that the invariant density $\rho_0$ is bounded by a constant $C$ in a neighborhood of $\Phi_0^t(1,1)$.  For $\eps$ sufficiently small, the set $R_{\eps}(t, t+\eps^{\alpha})$ is contained in this neighborhood. Using the diffeomorphism property of $\Phi_0^t$ mentioned above, it is easy to check that there is a number $C'>0$ such that for small $\eps$ the Lebesgue measure of $R_{\eps}(t,t+\eps^{\alpha})$ is bounded by $C'\eps^{\alpha + \beta}$. Therefore, $\mu_0(R_{\eps}(t,t+\eps^{\alpha}))\le CC'\eps^{\alpha + \beta}$ for $\eps$ small. 

Now we derive a lower bound for $\mu_0(R_{\eps}(t, t+\eps^{\alpha}))$. For small $\eps>0$, $\phi^{(1)}(\eps^{\alpha})< 1- \tfrac{\alpha}{2} \eps^{\alpha}$. Therefore, $R_{\eps}(0,\eps^{\alpha})\supset ((1- \tfrac{\alpha}{2} \eps^{\alpha},1)\times(1-\eps^{\beta}, 1))\cap\Gamma^{\circ}$. As $\mu_1$ is supported on $\Gamma$, this yields  
\begin{equation}    \label{eq:mu_1_lower_bd}
\mu_1(R_{\eps}(0,\eps^{\alpha})) \geq \int_{1- \frac{\alpha}{2} \eps^{\alpha}}^{\infty} \int_{1-\eps^{\beta}}^{\infty} \rho_1(y_1, y_2) \ d y_2 \ d y_1 \geq G_1(+\infty, +\infty) \frac{\lambda_0}{\lambda_1} \eps^{\lambda_1}, 
\end{equation}
where the second inequality follows from the proof of part (1), with the roles of $\lambda_0$ and $\lambda_1$ reversed.  Using~\eqref{eq:pdf},~\eqref{eq:R_identity} and the observation that $\Phi_0^{-s}: R_{\eps}(t, t+\eps^{\alpha}) \to R_{\eps}(t-s, t-s +\eps^{\alpha})$ is a diffeomorphism, we have  
\begin{align*}
 \mu_0(R_{\eps}(t,t+\eps^{\alpha}))
 &=\int_0^\infty \lambda_1 e^{-\lambda_0 s} \mu_1(R_{\eps}(t-s,t-s+\eps^{\alpha})) \ ds\\
 &\ge \int_0^\infty \lambda_1 e^{-\lambda_0 s} \mu_1(R_{\eps}(t-s,t-s+\eps^{\alpha})\cap R_{\eps}(0,\eps^{\alpha})) \ ds\\
 &= \int_0^\infty \lambda_1 e^{-\lambda_0 s} \mu_1(R_{\eps}((t-s,t-s+\eps^{\alpha})\cap(0,\eps^{\alpha}))) \ ds\\
 &\ge \int_{t-\eps^{\alpha}}^{t+\eps^{\alpha}} \lambda_1 e^{-\lambda_0 s} \mu_1(R_{\eps}((t-s,t-s+\eps^{\alpha})\cap(0,\eps^{\alpha})) \ ds \\
 &\ge c \int_{t-\eps^{\alpha}}^{t+\eps^{\alpha}} \mu_1(R_{\eps}((t-s,t-s+\eps^{\alpha})\cap(0,\eps^{\alpha}))) \ ds =: c A(\epsilon).  
\end{align*}
Here, $c > 0$ is a constant that doesn't depend on $\eps$.  To complete the estimate of $\mu_0(R_{\eps}(t,t+\eps^{\alpha}))$ from below, we use the lower bound on $\mu_1(R_{\eps}(0,\eps^{\alpha}))$ derived in~\eqref{eq:mu_1_lower_bd}. We have  
\begin{align*}
A(\eps) &=\int_{t-\eps^{\alpha}}^{t+\eps^{\alpha}} \int_{\R^2} \ONE_{\{(x_1,x_2)\in R_{\eps}((t-s,t-s+\eps^{\alpha})\cap(0, \eps^{\alpha}))\}} \ \mu_1(dx_1, d x_2) \ ds  \\
&=\int_{t-\eps^{\alpha}}^{t+\eps^{\alpha}} \int_{\R^2}\ONE_{\{x_1\in\phi^{(1)}((t-s,t-s+\eps^{\alpha})\cap(0,\eps^{\alpha}))\}}\ONE_{\{x_2\in I_{\eps}(x_1)\}} \ \mu_1(dx_1, dx_2) \ ds \\
&=\int_{\R^2} \ONE_{\{t(x_1) \in(0,\eps^{\alpha})\}} \ONE_{\{x_2\in I_{\eps}(x_1)\}} \int_{t-\eps^{\alpha}}^{t+\eps^{\alpha}} \ONE_{\{t(x_1)\in (t-s,t-s+\eps^{\alpha})\}} \ ds \ \mu_1(dx_1, dx_2) \\
&=\eps^{\alpha} \int_{\R^2} \ONE_{\{t(x_1)\in(0,\eps^{\alpha})\}} \ONE_{\{x_2\in I_{\eps}(x_1)\}} \ \mu_1(dx_1, dx_2) \\
&=\eps^{\alpha} \mu_1(R_{\eps}(0,\eps^{\alpha}))\ge G_1(+\infty, +\infty) \frac{\lambda_0}{\lambda_1} \eps^{\alpha+\lambda_1}.
\end{align*}
In the next to last line we used that $t(x_1)\in (t-s,t-s+\eps^{\alpha})$ if and only if $s\in (t-t(x_1),t-t(x_1)+\eps^{\alpha})$, and if $t(x_1) \in(0,\eps^{\alpha})$, one has $(t-t(x_1),t-t(x_1)+\eps^{\alpha}) \subset (t-\eps^{\alpha},t+\eps^{\alpha})$. 

Combining the resulting lower bound for $\mu_0(R_{\eps}(t, t + \eps^{\alpha}))$ with the upper bound derived earlier, we see that for sufficiently small $\eps$,
\[
c G_1(+\infty, +\infty) \frac{\lambda_0}{\lambda_1} \eps^{\alpha+\lambda_1}\le CC'\eps^{\alpha + \beta}.
\]
This is possible only if $\lambda_1 \ge \beta$. 
\epf

\bigskip

\section{Change of variables}        \label{sec:change_of_variables}

The fixed-point equations in \eqref{eq:Q_fixed_point} give integral
equations describing the invariant densities $\rho_i(x)$. These
expressions were 
obtained by applying \eqref{eq:pdf} multiple times. This captures the effect of pushing
forward the density until the time of the $n$th switch.  All of
these expressions are written as integral operators where the integrals
are taken over the length of the first $n$ exponential times.

The goal in this section is to take the expressions for the pushforwards
of the invariant densities in the case $n=2$ as integrals over the switching times and
change variables so that they can be viewed as integral operators over
the state space $\Gamma$. The precise version is given in
Lemma~\ref{lm:change_of_variables} at the end of this section.

It will be important to treat the points above, below and on the diagonal $\{x_1=x_2\}$
differently. To this end, let us introduce the sets 
\begin{align*}
\Gamma_l =& \{x \in \Gamma^{\circ}: x_1 < x_2\},  & \Gamma_r =& \{x \in \Gamma^{\circ}: x_1 > x_2\}, \\
\Gamma_m =& \{x \in \Gamma^{\circ}: x_1 = x_2\},  & \Gamma_s =& \Gamma_l \cup \Gamma_r. 
\end{align*}Recalling the definition of $T_i^n$ from \eqref{eq:TDef}, we define,
for $i \in \{0,1\}$ and $x \in \Gamma^{\circ}$, the following sets of
times
when the first switch, going backwards in time, occurred respectively to the right and to the left of the diagonal: 
\begin{align*}
R_i(x)  &= \{(s,t) \in T_i^2(x): \Psi_i^t(x) \in \Gamma_r\}, \\
L_i(x) &= \{(s,t) \in T_i^2(x):  \Psi_i^t(x) \in \Gamma_l\} \,. 
\end{align*} 
Observe that $R_1(x) = L_0((1,1)-x)$ and $L_1(x) = R_0((1,1) - x)$.
%
The next lemma provides the regularity needed to perform the desired
change of variables.

\begin{lemma}   \label{lm:F_injective}
For any $x \in \Gamma^{\circ}$ and $i \in \{0,1\}$, $(s,t) \mapsto \Psi_i^{(s,t)}(x)$ is a diffeomorphism from $R_i(x)$ onto $\Psi_i^{R_i(x)}(x)$ and from $L_i(x)$ onto $\Psi_i^{L_i(x)}(x)$.  
\end{lemma}

\bpf It suffices to prove the statement for $i=0$ because 
\begin{equation*}
\Psi_1^{(s,t)}(x) = (1, 1) - \Psi_0^{(s,t)}((1,1) - x), \quad (s,t) \in T_1^2(x). 
\end{equation*}
We only show that $(s,t) \mapsto \Psi_0^{(s,t)}(x)$ 
is a
diffeomorphism on $R_0(x)$, as the proof for
$L_0(x)$ is almost identical. We begin by showing that $(s,t) \mapsto \Psi_0^{(s,t)}(x)$
is injective on $R_0(x)$. To obtain a contradiction, suppose that this is not the case. Then, there
exist two distinct vectors $(s_1, t_1), (s_2,t_2) \in R_0(x)$ such
that $\Psi_0^{(s_1,t_1)}(x) = \Psi_0^{(s_2, t_2)}(x) =:y$. 
This implies that the $u_0$ trajectory through $x$ and the $u_1$ trajectory through $y$ intersect in two distinct points $z^{(1)} :=
\Psi_0^{t_1}(x)$ and $z^{(2)} := \Psi_0^{t_2}(x)$ that both
lie in $\Gamma_r$. For any two points $x, y \in \Gamma^{\circ}$, the $u_0$ trajectory through $x$ and the $u_1$ trajectory through $y$ intersect in at most two distinct points, so $z^{(1)}$ and $z^{(2)}$ are the only points
of intersection. Since $z^{(1)}, z^{(2)} \in \Gamma_r$, Lemma~\ref{lm:uniqueness_ac} 
implies that $\det U(z^{(i)})> 0$ for $i \in \{0,1\}$.
As a result, one trajectory
crosses the other in the same direction at both points of intersection, which
is impossible.

It remains to show that $\det \nabla_{(s,t)} \Psi_0^{(s,t)}(x) \neq 0$ for $(s,t) \in R_0(x)$. From~\eqref{eq:phi_v}, we derive 
\begin{equation}            \label{eq:explicit_representation}
\Psi_0^{(s,t)}(x) = \begin{pmatrix}
               1 - e^{\alpha s} + e^{\alpha (s+t)} x_1 \\
               1 - e^{\beta s} + e^{\beta (s+t)} x_2
              \end{pmatrix},
\end{equation}
which yields  
\begin{equation*}
 \nabla_{(s,t)} \Psi_0^{(s,t)}(x) = \begin{pmatrix}
                                            -\alpha e^{\alpha s} + \alpha x_1 e^{\alpha (s+t)} & \alpha x_1 e^{\alpha(s+t)} \\
                                             - \beta e^{\beta s} + \beta x_2 e^{\beta (s+t)} & \beta x_2 e^{\beta (s+t)} 
                                               \end{pmatrix}
\end{equation*}
and 
\begin{equation}     \label{eq:det_representation}
\det \nabla_{(s,t)} \Psi_0^{(s,t)}(x) = \alpha \beta e^{(\alpha + \beta) s} (x_1 e^{\alpha t} - x_2 e^{\beta t}) > 0.
\end{equation}
For the last inequality, we used that $\Psi_0^t(x) \in \Gamma_r$. 
\epf 

\bigskip
Let $i \in \{0,1\}$ and $x \in \Gamma^{\circ}$. We denote the inverse
of $(s,t) \mapsto \Psi_i^{(s,t)}(x)$ as a map from $R_i(x)$ onto $\Psi_i^{R_i(x)}(x)$ by $\chi^{r,x}_i$, and the inverse of $(s,t)
\mapsto\Psi_i^{(s,t)}(x)$ as a map from  $L_i(x)$ onto $\Psi_i^{L_i(x)}(x)$ by $\chi^{l,x}_i$. 
With $\lambda_i(2)$ and $\lambda_i^{(2)}$ defined as in~\eqref{eq:lambda_def_1} and~\eqref{eq:lambda_def_2}, respectively, we also introduce the functions 
\begin{align*}
f_i(\tbf,x) =& \lambda_i(2) e^{- \langle \lambda_i^{(2)}, \tbf \rangle} J_i^{\tbf}(x), \quad \tbf \in R_i(x) \cup L_i(x), \\
K^r_i(x, y) =& f_i(\chi_i^{r,x}(y), x)  \lvert \det \nabla_y \chi_i^{r,x}(y) \rvert, \quad y \in \Psi_i^{R_i(x)}(x), \\
K^l_i(x,y) =& f_i(\chi_i^{l,x}(y),x)  \lvert \det \nabla_y \chi_i^{l,x}(y) \rvert, \quad y \in \Psi_i^{L_i(x)}(x).  
\end{align*} 

We are now in a position to perform the desired change of variables on
the operator $\Qc_2$ which was defined in \eqref{eq:QcDef} and used in the fix
point equations \eqref{eq:Q_fixed_point}. 
\begin{lemma}    \label{lm:change_of_variables}
For any $i \in \{0,1\}$ and $x \in \Gamma^{\circ}$,  
\begin{equation}    \label{eq:change_of_variables}
\int_{R_i(x)} \lambda_i(2) e^{-\langle \lambda_i^{(2)}, \tbf \rangle} J^{\tbf}_i(x) \rho_i(\Psi^{\tbf}_i x) \ d \tbf = \int_{\Psi_i^{R_i(x)}(x)} \rho_i(y) K^r_i(x,y)  \ dy.  
\end{equation}
In~\eqref{eq:change_of_variables}, one can replace $R_i(x)$ and $K^r_i$ together with $L_i(x)$ and $K^l_i$.  
\end{lemma}

\bpf The formula follows after applying the change of variables
$y = \Psi_i^{\tbf}(x)$ justified by Lemma~\ref{lm:F_injective}. 
\epf 

\bigskip

\section{Boundedness}        \label{sec:proof_bounded}

In this section, we prove Theorem~\ref{thm:bounded} that describes under which conditions the invariant density $\rho_0$ stays bounded. The main difficulties in proving this result stem from two sources: the exponential contraction in the vicinity of the critical points $(0,0)$ and $(1,1)$, and the fact, exhibited in Lemma~\ref{lm:uniqueness_ac}, that the vector fields $u_0$ and $u_1$ are collinear at every point on the line $x_1 = x_2$. As we saw in Section~\ref{sec:proof_unbounded}, exponential contraction is an essential problem that gives rise to singularities of the invariant densities for slow switching. The lack of ellipticity along the diagonal $x_1 = x_2$ creates technical challenges because switches close to the diagonal have a less pronounced regularizing effect on the invariant densities.  At the same time, switches close to the diagonal do not actively spoil the densities as long as they occur sufficiently far from the two critical points.  

Throughout this section, we will use the following basic facts about the switching system, at times without explicitly referring to them. 

\begin{lemma}         \label{lm:basics}
For any $x \in \Gamma^{\circ}$, the following statements hold. 
\begin{enumerate}
\item For any $i \in \{0,1\}$ there is a unique $\theta_i(x) \in \R$ such that $\det U(\Psi_i^{\theta_i(x)} x) = 0$. We have $\Psi_i^{\theta_i(x)}(x) \in \Gamma^{\circ}$. 
\item We have 
$$
\frac{d}{dt} \det U(\Psi_i^t x) = \alpha \beta \left((i-x_2) \beta e^{\beta t} - (i-x_1) \alpha e^{\alpha t} \right), \quad i \in \{0,1\}, 
$$
and
\begin{equation*}
\frac{d}{dt} \det U \left( \Psi_0^t x \right) \vert_{t = 0} > \beta \det U(x), \quad \frac{d}{dt} \det U\left(\Psi_1^t x \right) \vert_{t=0} < \alpha \det U(x). 
\end{equation*}
In particular, if $x_1 = x_2$ and if $\eps > 0$, there is a unique $t_x(\eps) > 0$ such that $\det U(\Psi_0^{t_x(\eps)} x) = \eps$. 
\item Suppose now that $x_1 = x_2$, and let $y \in \Gamma^{\circ}$ such that $x_1 < y_1 = y_2$. Then, 
$$ 
[\Psi_0^{t_y(\eps)}(y)]_2 - [\Psi_0^{t_x(\eps)}(x)]_2 > 0, \quad \eps > 0.
$$ 
Here, $[z]_2 := z_2$ for $z = (z_1, z_2) \in \R^2$. 
\end{enumerate}
\end{lemma}

We omit the proof of this lemma.   To illustrate our main strategy for establishing Theorem~\ref{thm:bounded}, we first show that $\rho_0$ is bounded on the part of $\Gamma$ that lies below the diagonal $x_1 = x_2$. This statement has a comparatively simple proof as we do not need to address the lack of ellipticity and as there is no danger of entering regions with strong exponential contraction. 

\begin{proposition}     \label{prop:below_diag} 
Let $K \subset \Gamma_r \cup \partial \Gamma_1$ be compact. Then, $\rho_0$ is bounded on $K$ for every $\lambda_0, \lambda_1 > 0$. 
\end{proposition}

\bpf  Using~\eqref{eq:Q_fixed_point} for $n = 2$, we
have  
\begin{equation}     \label{eq:integral_equation_2}
 \rho_0(x) = \int_{R_0(x)} \lambda_0(2)  e^{-\langle \lambda_0^{(2)}, \tbf \rangle}
J_0^{\tbf}(x) \rho_0(\Psi_0^{\tbf}x) \ d \tbf, \quad x \in \Gamma_r.
\end{equation}
Lemma~\ref{lm:change_of_variables} then yields 
\begin{equation}    \label{eq:bounded}
 \rho_0(x) = \int_{\Psi_0^{R_0(x)}(x)} \rho_0(y) K^r_0(x, y) \ dy. 
\end{equation}
Since $\rho_0$ is integrable, it is enough to show that there is $c > 0$ such that 
$$ 
K_0^r(x,y) \leq c, \quad x \in K \cap \Gamma_r, \ y \in \Psi_0^{R_0(x)}(x). 
$$ 
Let $x \in K \cap \Gamma_r$, $y \in \Psi_0^{R_0(x)}(x)$, and set $\tbf = (s,t) = \chi_0^{r,x}(y)$. Since $\Psi_0^t(x) \in \Gamma_r$, the matrix $U(\Psi_0^t x)$ is invertible. In~\cite{torus}, proof of Theorem 2, the formula 
$$ 
\nabla_x \Psi^{(s,t)}_0(x) = -\nabla_{(s,t)} \Psi^{(s,t)}_0(x) U(\Psi_0^t x)^{-1} \nabla_x \Psi_0^t(x)  
$$ 
was established. It shows how the effect of variations in the initial point on the
final point after two switches (the lefthand side) can be translated
into an equivalent variation in the switching times, with the translation given by the two rightmost
terms on the righthand side. 
If we set $z = \Psi_0^t(x)$, the formula above yields 
$$
J_0^{(s,t)}(x) = \det \nabla_{(s,t)} \Psi_0^{(s,t)}(x) \det U(z)^{-1} \det \nabla_x \Psi_0^t(x). 
$$
Hence, since $\det \nabla_y \chi_0^{r,x}(y) = (\det \nabla_{(s,t)} \psi_0^{(s,t)}(x))^{-1}$, 
$$
K_0^r(x,y) = \lambda_0(2) e^{-\langle \lambda_0^{(2)}, \tbf \rangle} \det U(z)^{-1} \det \nabla_x \Psi_0^t(x) = \lambda_0(2)  e^{-\lambda_1 s} \frac{e^{(\alpha + \beta -\lambda_0)t}}{\alpha \beta (z_1 - z_2)}. 
$$
By definition, $t$ is the time it takes to move backward along
the $u_0$ trajectory from $x$ to $z$. As $z \in \Gamma^{\circ}$ and thus $z_2 < 1$, we have  
\begin{equation*}
t < -\frac{1}{\beta} \ln(x_2)  \leq \sup_{x \in K} \left(-\frac{1}{\beta} \ln(x_2) \right) < \infty, 
\end{equation*}
where one should note that $K$ is compact and only contains points $x = (x_1, x_2)$ such that $x_2 > 0$. 
As $z$ lies on the backward
$u_0$ trajectory starting at $x$, we also have 
\begin{equation*}
\frac{1}{z_1 - z_2} \leq \frac{1}{x_1 - x_2} \leq \sup_{x \in K} \frac{1}{x_1 - x_2} < \infty.
\end{equation*}
This completes the proof. 
\epf

\bigskip

\subsection{Switches close to $\Gamma_m$}     \label{ssec:diagonal}

If $y \in \Psi_0^{R_0(x)}(x)$ is chosen in such a way that $z$, the point where the switch from $u_1$ to $u_0$ occurs (cf. proof of Proposition~\ref{prop:below_diag}), is close to the diagonal $\Gamma_m$, the term $(z_1 - z_2)^{-1}$ is very large, and so is $K^r_0(x, y)$. For $x \in \Gamma_l \cup \Gamma_m$, the point $z$ can become arbitrarily close to $\Gamma_m$, which prevents $K^r_0(x, \cdot)$ from being bounded on $\Psi_0^{R_0(x)}(x)$.  The proof of Proposition~\ref{prop:below_diag} can therefore not be extended to the case of $x \in \Gamma_l \cup \Gamma_m$ in a straightforward way. In this subsection, we describe an approach for dealing with this difficulty. 

For any $n \in \N$, $x \in \Gamma^{\circ}$, and $\epsilon > 0$, let 
\begin{align*}
M^{\epsilon}_n(x) =& \left\{ \tbf \in T_0^n(x):  \left \lvert \det U(\Psi_0^{(t_{n-j}, \ldots, t_n)} x) \right \rvert < \epsilon, 0 \leq j \leq n-1 \right\}, \\
S^{\epsilon}_n(x)=& \left\{\tbf \in T_0^{n+1}(x): (t_3, \ldots, t_{n+1}) \in M^{\epsilon}_{n-1}(x), \left \lvert \det U(\Psi_0^{(t_2, \ldots, t_{n+1})} x) \right \rvert > \epsilon \right\}.
\end{align*} 
The condition $(t_3, \ldots, t_{n+1}) \in M_{n-1}^{\epsilon}(x)$ is void for $n=1$. For $h \in L^1(\Gamma^{\circ})$, we define 
\begin{align*}
\AC^{\epsilon}_n h(x) =& \int_{S^{\epsilon}_n(x)} \lambda_0(n+1) e^{-\langle \lambda_0^{(n+1)}, \tbf \rangle} J_0^{\tbf}(x) h(\Psi_0^{\tbf} x) \ d \tbf, \\
\Bc^{\epsilon}_n h(x) =&  \int_{M^{\epsilon}_n(x)} \lambda_0(n) e^{-\langle \lambda_0^{(n)}, \tbf \rangle} J_0^{\tbf}(x) h(\Psi_0^{\tbf} x) \ d \tbf. 
\end{align*}
To avoid distinguishing between the cases of even $n$ and odd $n$, we introduce the shorthand 
$$ 
i_n = \begin{cases}
            0, & \quad n \equiv 0 \mod 2, \\
            1, & \quad n \equiv 1 \mod 2. 
\end{cases}
$$

\begin{lemma}     \label{lm:Hennion}
 For $n \in \N$, $x \in \Gamma^{\circ}$, and $\epsilon > 0$, 
\begin{equation*}
 \rho_0(x) = \sum_{k=1}^n \AC^{\epsilon}_k \rho_{i_{k+1}}(x) + \Bc^{\epsilon}_n \rho_{i_n}(x).  
\end{equation*}
\end{lemma}

\bpf The proof is by induction. The formula in~\eqref{eq:Q_fixed_point} gives 
\begin{equation}   \label{eq:integral_formula_one_switch}
 \rho_0(x) = \int_{T_0^1(x)} \lambda_1 e^{-\lambda_0 t} J_0^t(x) \rho_1(\Psi_0^t x) \ dt. 
\end{equation}
Statement (2) in Lemma~\ref{lm:basics} implies that for $y \in \R^2$, $t \in \R$, and $i \in \{0,1\}$, 
$$
\frac{d}{dt} \det U(\Psi_i^t y) = 0 
$$
if and only if 
$$
t = \frac{1}{\beta - \alpha} \ln \left(\frac{\alpha (y_1 - i)}{\beta (y_2 - i)} \right).
$$ 
This shows that the set $T_0^1(x)$ is, up to a set of Lebesgue measure zero, the disjoint union of $M_1^{\epsilon}(x)$ and the set of $t \in \R_+$ such that $(s,t) \in S_1^{\epsilon}(x)$ for some $s \in \R_+$. Thus, the right side of~\eqref{eq:integral_formula_one_switch} can be written as 
\begin{equation}    \label{eq:Hennion_2}
\int_{t: \exists s \ \text{s.t.} \ (s,t) \in S_1^{\epsilon}(x)} \lambda_1 e^{-\lambda_0 t} J_0^t(x) \rho_1(\Psi_0^t x) \ d t + \Bc^{\epsilon}_1 \rho_1(x). 
\end{equation}
In complete analogy to~\eqref{eq:integral_formula_one_switch}, we have 
\begin{equation}     \label{eq:Hennion_3} 
\rho_1(y) = \int_{T_1^1(y)} \lambda_0 e^{-\lambda_1 s} J_1^s(y) \rho_0(\Psi^s_1 y) \ ds. 
\end{equation} 
If we plug this identity into the first summand in~\eqref{eq:Hennion_2}, we obtain the desired formula in the base case $n = 1$. 

In the induction step, assume the formula holds for some $n \in \N$. With the notation $\tbf = (t_2, \ldots, t_{n+1})$, we can write 
\begin{align} \label{eq:Hennion_4}
&\Bc^{\epsilon}_n \rho_{i_n}(x) \\
=& \int_{M_n^{\epsilon}(x)}  \int_{T^1_{i_n}(\Psi^{\tbf}_0 x)}  \lambda_0(n+1) e^{- \langle \lambda_0^{(n+1)}, (t_1, \tbf)^{\top} \rangle} J_0^{(t_1, \tbf)}(x) \rho_{i_{n+1}}(\Psi_0^{(t_1, \tbf)} x) \ dt_1 \ d \tbf. \notag
\end{align}
The set
\begin{equation*}
\left \{(t_1, \tbf) \in \R^{n+1}_{+}: \tbf \in M_n^{\epsilon}(x), \ t_1 \in T^1_{i_n}(\Psi_0^{\tbf} x)\right\}
\end{equation*}
is, again up to a set of Lebesgue measure zero, the disjoint union of $M_{n+1}^{\epsilon}(x)$ and the set of $(t_1, \tbf) \in \R_{+}^{n+1}$ such that $(t_0, t_1, \tbf) \in S_{n+1}^{\epsilon}(x)$ for some $t_0 \in \R_{+}$. Therefore, the right side of~\eqref{eq:Hennion_4} becomes 
$$
\Bc_{n+1}^{\epsilon} \rho_{i_{n+1}}(x) + \int_{\tbf: \exists t_0 \ \text{s.t.} \ (t_0, \tbf) \in S_{n+1}^{\epsilon}(x)} \lambda_0(n+1) e^{-\langle \lambda_0^{(n+1)}, \tbf \rangle} J_0^{\tbf}(x) \rho_{i_{n+1}}(\Psi_0^{\tbf} x) \ d \tbf,  
$$
where we have set $\tbf = (t_1, \ldots, t_{n+1})$. It remains to show
that the integral term above equals $\AC_{n+1}^{\epsilon} \rho_{i_{n+2}}$. 
This follows from plugging~\eqref{eq:integral_formula_one_switch} or~\eqref{eq:Hennion_3} into said integral term.  
\epf

\bigskip

Next, we show that for large $n$, the contribution of $\Bc_n^{\epsilon} \rho_{i_n}$ in the formula from Lemma~\ref{lm:Hennion} is small.

\begin{lemma}       \label{lm:convergence_to_0}
We have $\lim_{n \to \infty} \|\Bc^{\epsilon}_n \|_{op} = 0$, where $\|\cdot\|_{op}$ is the operator norm for operators on $L^1(\Gamma^{\circ})$. 
\end{lemma}

\bpf  For fixed $n \in \N$ and $h \in L^1(\Gamma^{\circ})$, we have 
\begin{equation}    \label{eq:con_0_1}
\| \Bc^{\epsilon}_n h \|_{L^1} \leq \int_{\widehat{M^{\epsilon}_n}} \lambda_0(n) e^{- \langle \lambda_0^{(n)}, \tbf \rangle} J_0^{\tbf}(x) \lvert h(\Psi_0^{\tbf} x) \rvert \ dx \ d \tbf, 
\end{equation}
where 
$$ 
\widehat{M^{\epsilon}_n} = \{(x,\tbf) \in \Gamma^{\circ} \times \R^n_+: \tbf \in M^{\epsilon}_n(x)\}. 
$$ 
Let us set 
\begin{equation*}
M_n^{\epsilon} = \bigcup_{x \in \Gamma^{\circ}} M_n^{\epsilon}(x). 
\end{equation*}
Since $\widehat{M_n^{\epsilon}} \subset \Gamma^{\circ} \times M_n^{\epsilon}$, the right side of~\eqref{eq:con_0_1} is bounded by 
\begin{equation*}
\int_{M_n^{\epsilon}} \lambda_0(n) e^{- \langle \lambda_0^{(n)}, \tbf \rangle} \int_{\Gamma^{\circ}} J_0^{\tbf}(x) \lvert h(\Psi_0^{\tbf} x) \rvert \ dx \ d\tbf.
\end{equation*}
With the change of variables $y = \Psi_0^{\tbf}(x)$, the expression above becomes  
\begin{equation*}
\|h \|_{L^1} \int_{M_n^{\epsilon}} \lambda_0(n) e^{-\langle \lambda_0^{(n)}, \tbf \rangle} \ d \tbf.  
\end{equation*}
Thus, 
$$ 
\|\Bc_n^{\epsilon}\|_{op} \leq \int_{M_n^{\epsilon}} \lambda_0(n) e^{-\langle \lambda_0^{(n)}, \tbf \rangle} \ d \tbf =: b_n. 
$$ 
As $M_{n+1}^{\epsilon} \subset \R_{+} \times M_n^{\epsilon}$, we have 
\begin{equation*}
b_{n+1} \leq \int_{\R_{+}} \lambda_{i_{n+1}}  e^{-\lambda_{i_n} t} \ dt \int_{M_n^{\epsilon}} \lambda_0(n) e^{-\langle \lambda_0^{(n)}, \tbf \rangle} \ d \tbf = \frac{\lambda_{i_{n+1}}}{\lambda_{i_n}} b_n,
\end{equation*}
so $\tfrac{b_{n+1}}{b_n}$ is bounded. To show that $\lim_{n \to \infty} b_n = 0$, it then suffices to show that there are $c \in (0,1)$ and $N \in \N$ such that 
$$
b_{2n + 3} \leq c b_{2n + 1}, \quad n \geq N. 
$$
We claim that there are $\tau_1, \tau_2, T > 0$ such that for $n$ sufficiently large, 
\begin{equation}   \label{eq:big_set}
M_{2n+3}^{\epsilon} \subset \left( \left( \R_{+} \times \left((0, \tau_1] \cup [\tau_2, \infty) \right) \right) \cup \left( (0, T) \times \left( \tau_1, \tau_2 \right) \right) \right) \times M_{2n + 1}^{\epsilon}.
\end{equation}
The idea behind this claim is the following:  If the time between two consecutive switches that both happen close to the diagonal $\Gamma_m$ is neither very short nor very long, i.e. if it falls within $(\tau_1, \tau_2)$, then, at the time of the second switch, the switching trajectory cannot end up close to the critical point of the vector field to which the second switch is made. As a result, the time spent in $\Gamma^{\circ}$ after the second switch, following the time-reversed flow, is bounded by $T$.   To prove the claim in~\eqref{eq:big_set}, we first notice that for all $\tau > 0$  
$$ 
\det U(\Phi_1^{\tau}(0,0)) \neq 0. 
$$ 
Let us pick one such $\tau$. Since the map 
$$ 
(x,t) \mapsto \det U(\Phi_1^t x)
$$ 
is jointly continuous, there are $\hat \epsilon > 0$, $\tau_2 > \tau_1 > 0$ and a neighborhood $V$ of $(0,0)$ such that 
$$ 
\lvert \det U(\Phi_1^t x) \rvert > \hat \epsilon, \quad t \in (\tau_1, \tau_2), \ x \in V. 
$$ 
Then, set 
$$ 
T = 1 + \sup\{t \geq 0: \exists x \in V^c \cap \Gamma^{\circ} \ \text{s.t.} \ \Psi_0^t (x) \in \Gamma^{\circ}\}. 
$$ 
This defines a finite quantity because $x_2 e^{\beta t} > 1$ for $t > - \tfrac{1}{\beta} \ln(x_2)$, and \\ $\sup_{x \in V^c \cap \Gamma^{\circ}} -\ln(x_2) < \infty$. Let $(t_1, \ldots, t_{2n+3}) \in M_{2n+3}^{\epsilon}$. Then, $(t_1, \ldots, t_{2n+3}) \in M_{2n+3}^{\epsilon}(x)$ for some $x \in \Gamma^{\circ}$, and in particular $(t_3, \ldots, t_{2n+3}) \in M_{2n+1}^{\epsilon}(x)$. Suppose that $t_2 \in (\tau_1, \tau_2)$ and set 
$$ 
z = \Psi^{(t_2, \ldots, t_{2n+3})}(x). 
$$ 
We claim that $z \in V^c$. If $z$ was an element of $V$, we would have 
$$ 
\hat \epsilon < \lvert \det U(\Phi_1^{t_2} z) \rvert = \lvert \det U(\Psi_0^{(t_3, \ldots, t_{2n+3})} x) \rvert. 
$$ 
On the other hand, as $(t_3, \ldots, t_{2n+3}) \in M_{2n+1}^{\epsilon}(x)$, we have 
$$ 
\lvert \det U(\Psi_0^{(t_3, \ldots, t_{2n+3})} x) \rvert < \epsilon \leq \hat \epsilon 
$$ 
for $\epsilon$ sufficiently small, a contradiction.  

From $z \in V^c$ it follows that $\Psi_0^t(z) \notin \Gamma^{\circ}$ for $t \geq T$. Hence, $t_1 < T$, and the inclusion in~\eqref{eq:big_set} is proved. As a result,   
$$
b_{2n + 3} \leq c b_{2n+1}, 
$$
where 
\begin{equation*}
c = \int_{(0, \tau_1] \cup [\tau_2, \infty)} \lambda_1 e^{-\lambda_1 t} \ dt + \int_0^T \lambda_0 e^{-\lambda_0 s} \ ds \int_{\tau_1}^{\tau_2} \lambda_1 e^{-\lambda_1 t} \ dt < 1.
\end{equation*}
This completes the proof of the lemma.  
\epf

\bigskip

\begin{corollary}    \label{co:convergence_to_0}
For $\epsilon > 0$ sufficiently small, we have $\lim_{n \to \infty} \Bc_n^{\epsilon} \rho_{i_n}(x) = 0$ for Lebesgue almost every $x \in \Gamma^{\circ}$. 
\end{corollary}

\bpf Recall from the proof of Lemma~\ref{lm:Hennion} that 
\begin{align*}
&\Bc_n^{\epsilon} \rho_{i_n}(x) \\
=& \Bc_{n+1}^{\epsilon} \rho_{i_{n+1}}(x) + \int_{\tbf: \exists t_0 \ \text{s.t.} \ (t_0,\tbf) \in S_{n+1}^{\epsilon}(x)} \lambda_0(n+1) e^{- \langle \lambda_0^{(n+1)}, \tbf \rangle} J_0^{\tbf}(x) \rho_{i_{n+1}}(\Psi_0^{\tbf} x) \ d \tbf. 
\end{align*} 
Thus, $(\Bc_n^{\epsilon} \rho_{i_n})_{n \geq 1}$ is a pointwise monotone decreasing sequence of nonnegative functions. To prove the corollary, it is therefore enough to show that $(\Bc_n^{\epsilon} \rho_{i_n})_{n \geq 1}$ converges to $0$ in $L^1(\Gamma^{\circ})$, which follows immediately from Lemma~\ref{lm:convergence_to_0}. 
\epf 

\bigskip

The next lemma is a counterpart to the change of variables - formula in Lemma~\ref{lm:change_of_variables}. 
For $i \in \{0,1\}$, $x \in \Gamma^{\circ}$, and $\epsilon > 0$, let 
\begin{align*}
R_i^{\eps}(x) =& \{(s,t) \in T_i^2(x): \det U(\Psi_i^t x) > \eps\}, \\
L_i^{\eps}(x) =& \{(s,t) \in T_i^2(x): \det U(\Psi_i^t x) < -\eps\}. 
\end{align*}
In addition, we define 
\begin{align*}
\Ic^{r,\eps}_i(x) &=  \int_{\Psi_i^{R_i^{\eps}(x)}(x)} \rho_i(y) K^r_i(x, y) \ dy, \\
\Ic^{l,\eps}_i(x) &= \int_{\Psi_i^{L_i^{\eps}(x)}(x)} \rho_i(y) K^l_i(x, y) \ dy, \\
\Ic_i^{\epsilon}(x) &= \Ic^{r,\eps}_i(x) + \Ic^{l,\eps}_i(x).  
\end{align*}

\begin{lemma}     \label{lm:change_of_variables_1}
Let $x \in \Gamma^{\circ}$ and $\epsilon > 0$. Then,  
\begin{equation*}  
\AC_1^{\eps} \rho_{i_2}(x) = \Ic_0^{\eps}(x). 
\end{equation*}
For $k > 1$, we have 
\begin{equation}    \label{eq:cov_1}
\AC_k^{\eps} \rho_{i_{k+1}}(x) = \int_{M^{\eps}_{k-1}(x)} \lambda_0(k-1) e^{ \langle (\alpha + \beta) \mathbbm{1} - \lambda_0^{(k-1)}, \tbf \rangle} \Ic_{i_{k+1}}^{\epsilon}(\Psi_0^{\tbf} x) \ d \tbf, 
\end{equation} 
where $\mathbbm{1} = (1,1, \ldots, 1)^{\top}$. 
\end{lemma}

\bpf Observe that
\begin{align}       
\AC_1^{\eps} \rho_{i_2}(x) =& \int_{R_0^{\eps}(x)} \lambda_0(2) e^{-\langle \lambda_0^{(2)}, \tbf \rangle} J_0^{\tbf}(x) \rho_{i_2}(\Psi_0^{\tbf} x) \ d \tbf  \label{eq:cov_2} \\
& + \int_{L_0^{\eps}(x)} \lambda_0(2) e^{- \langle \lambda_0^{(2)}, \tbf \rangle} J_0^{\tbf}(x) \rho_{i_2}(\Psi_0^{\tbf} x) \ d\tbf, \notag   
\end{align}
and that Lemma~\ref{lm:change_of_variables} continues to hold if one replaces $R_0(x)$ with $R_0^{\eps}(x)$ and $L_0(x)$ with $L_0^{\eps}(x)$. The change of variables--formula in Lemma~\ref{lm:change_of_variables}
then yields that the right 
side of~\eqref{eq:cov_2} equals $\Ic_0^{\eps}(x)$. 

For $k > 1$, write $\AC_k^{\eps} \rho_{i_{k+1}}(x)$ as 
\begin{align}     \label{eq:cov_4}
&  \int_{M_{k-1}^{\epsilon}(x)} \ d t_3 \ldots d t_{k+1} \ \lambda_0(k-1) e^{- \langle \lambda_0^{(k-1)}, (t_3, \ldots, t_{k+1})^{\top} \rangle} \\
& \biggl( \int_{R_{i_{k+1}}^{\eps}(\Psi_0^{(t_3, \ldots, t_{k+1})} x)} \ d t_1 \ d t_2 \ \lambda_{i_{k+1}}(2) e^{- \langle \lambda_{i_{k+1}}^{(2)}, (t_1, t_2)^{\top} \rangle}J_0^{\tbf}(x) \rho_{i_{k+1}}(\Psi_0^{\tbf} x) \notag \\
& + \int_{L_{i_{k+1}}^{\eps}(\Psi_0^{(t_3, \ldots, t_{k+1})} x)} \ d t_1 \ d t_2 \ \lambda_{i_{k+1}}(2) e^{- \langle \lambda_{i_{k+1}}^{(2)}, (t_1, t_2)^{\top} \rangle} J_0^{\tbf}(x) \rho_{i_{k+1}}(\Psi_0^{\tbf} x) \biggr),    \notag
\end{align}
where $\tbf = (t_1, \ldots, t_{k+1})^{\top}$. 
Since 
$$ 
J_0^{\tbf}(x) = J_{i_{k+1}}^{(t_1, t_2)}(\Psi_0^{(t_3, \ldots, t_{k+1})} x) J_0^{(t_3, \ldots, t_{k+1})}(x) = J_{i_{k+1}}^{(t_1, t_2)}(\Psi_0^{(t_3, \ldots, t_{k+1})} x) e^{(\alpha + \beta) \sum_{j=3}^{k+1} t_j}, 
$$ 
we obtain the desired formula after applying Lemma~\ref{lm:change_of_variables} to the integrals in the second and third line of~\eqref{eq:cov_4}.
\epf 

\bigskip

\begin{lemma}      \label{lm:estimates}
For $i \in \{0,1\}$, $\epsilon > 0$, and $x \in \Gamma^{\circ}$, we have
\begin{equation*}  
\Ic_i^{r,\eps}(x) \leq \begin{cases}
                                      \frac{\lambda_0 \lambda_1}{\epsilon}, & \quad \lambda_i \geq \alpha + \beta, \\
                                      \frac{\lambda_0 \lambda_1}{\eps} e^{(\alpha + \beta - \lambda_i) \tau_i(x)}, & \quad \lambda_i < \alpha + \beta, 
\end{cases} 
\end{equation*}
where 
$$ 
\tau_i(x) = \sup \{t \geq 0: \Psi_i^t(x) \in \Gamma^{\circ}\}. 
$$
The estimate continues to hold if one replaces $\Ic_i^{r,\eps}(x)$ with $\Ic_i^{l, \eps}(x)$.  
\end{lemma}

\bpf Setting $(s,t) = \chi_i^{r,x}(y)$ and $z = \Psi_i^t(x)$, we have 
$$ 
K_i^r(x,y) = \lambda_i(2) e^{- \langle \lambda_i^{(2)}, (s, t)^{\top} \rangle} \lvert \det U(z)^{-1} \rvert \det \nabla_x \Psi_i^t(x) = \lambda_i(2) e^{-\lambda_{1-i} s} \frac{e^{(\alpha + \beta-\lambda_i) t}}{\lvert \det U(z) \rvert}. 
$$ 
Since $\lvert \det U(z) \rvert \geq \eps$ for $y \in \Psi_i^{R_i^{\eps}(x)}(x)$, the desired estimate follows.  
\epf

\bigskip

For $n \in \N$, $\eps > 0$, $x \in \Gamma^{\circ}$, and $\lambda_0, \lambda_1 > 0$, set 
$$ 
\M_n^{\eps,\lambda_0, \lambda_1}(x) = \int_{M_n^{\eps}(x)} e^{- \langle \lambda_0^{(n)}, \tbf \rangle} \ d \tbf. 
$$ 

\begin{lemma}    \label{lm:geometric_series}
For any $\lambda_0, \lambda_1 > 0$, there is a function $f(\epsilon)$ such that $\lim_{\epsilon \downarrow 0} f(\epsilon) = 0$ and 
\begin{equation*}
\M_{n+2}^{\epsilon, \lambda_0, \lambda_1}(x) \leq f(\epsilon) \M_n^{\epsilon, \lambda_0, \lambda_1}(x), \quad n \in \N, \ x \in \Gamma^{\circ}. 
\end{equation*}
\end{lemma}

\bpf   For $\eps > 0$ and $i \in \{0,1\}$, let 
$$ 
\Gamma_i^{\eps} = \left\{x \in \Gamma^{\circ}: \lvert \det U(x) \rvert < \eps, \ (-1)^i x_1 > \frac{1}{2}-i \right\}. 
$$ 
For $n \in \N$, $x \in \Gamma^{\circ}$, and $\epsilon > 0$, we have 
$$
\M_{n+2}^{\epsilon, \lambda_0, \lambda_1}(x) =  \int_{M_n^{\eps}(x)}  e^{- \langle \lambda_0^{(n)}, \tbf \rangle} \sum_{(i,j) \in \{0,1\}^2} I_{i,j}^{\tbf} \ d \tbf, 
$$ 
where 
$$ 
I_{i,j}^{\tbf} = \int_{t_2 \in \R_{+}: \ \Psi_0^{(t_2, \tbf)}(x) \in \Gamma_j^{\epsilon}} \ d t_2 \ e^{- \lambda_{i_n} t_2} \int_{t_1 \in \R_{+}: \ \Psi_0^{(t_1, t_2, \tbf)}(x) \in \Gamma_i^{\epsilon}}  \ d t_1 \ e^{- \lambda_{i_{n+1}} t_1}. 
$$
For $\tbf \in M_n^{\eps}(x)$ and $(i,j) \in \{0,1\}^2$, we now derive an upper bound on $I_{i,j}^{\tbf}$. For $x \in \Gamma_0^{\eps}$, $t \geq 0$, and $\eps < \tfrac{\alpha}{2} (\alpha - \beta)$, part (2) of Lemma~\ref{lm:basics} yields 
\begin{align*} 
\frac{d}{dt} \det U(\Psi_0^t x) <& e^{\beta t} \beta  (- \alpha (\alpha - \beta) x_1 - \det U(x)) \\
\leq& e^{\beta t} \beta  \left(-\frac{\alpha}{2} (\alpha - \beta) + \eps \right) \leq \beta \left(-\frac{\alpha}{2} (\alpha - \beta) + \eps \right).   
\end{align*}
Similarly, 
$$ 
\frac{d}{dt} \det U(\Psi_1^t x)  > \beta \left(\frac{\alpha}{2} (\alpha - \beta) - \eps \right), \quad x \in \Gamma_1^{\eps}, \ t \geq 0. 
$$ 
Thus, we have for $x \in \Gamma_i^{\eps}$ and 
$$
t \geq \frac{2 \eps}{\beta (\frac{\alpha}{2} (\alpha - \beta) -\eps)}=: \tau_{\eps}
$$
the estimate
\begin{align}    
\lvert \det U(\Psi_i^t x) \rvert \geq& \lvert \det U(\Psi_i^t x) - \det U(x) \rvert - \lvert \det U(x) \rvert \label{eq:exit_time_estimate} \\
>& \beta \left(\frac{\alpha}{2} (\alpha - \beta) - \eps \right) t - \eps \geq \eps.    \notag
\end{align} 
We distinguish between two cases. Suppose first that $j = i_{n+1}$.  With~\eqref{eq:exit_time_estimate}, we obtain for any $t_2 \in \R_{+}$ such that $\Psi_0^{(t_2, \tbf)}(x) \in \Gamma_j^{\eps}$ the estimate 
\begin{equation*}
\int_{t_1 \in \R_{+}: \Psi_0^{(t_1, t_2, \tbf)}(x) \in \Gamma_i^{\eps}} e^{- \lambda_{i_{n+1}} t_1} \ dt_1 \leq \int_0^{\tau_{\epsilon}} e^{- \lambda_{i_{n+1}} t_1} \ dt_1 \leq \tau_{\epsilon}. 
\end{equation*}
Then, 
$$
I_{i,j}^{\tbf} \leq \tau_{\epsilon} \int_0^{\infty} e^{- \lambda_{i_n} t_2} \ dt_2 = \frac{\tau_{\epsilon}}{\lambda_{i_n}}. 
$$
If $j = i_n$, we first use the obvious estimate 
\begin{equation}    \label{eq:case_2_est}
I_{i,j}^{\tbf} \leq \frac{1}{\lambda_{i_{n+1}}} \int_{t_2 \in \R_{+}: \Psi_0^{(t_2, \tbf)}(x) \in \Gamma_j^{\epsilon}} e^{- \lambda_{i_n} t_2} \ dt_2. 
\end{equation}
If $\{t_2 \in \R_{+}: \Psi_0^{(t_2, \tbf)}(x) \in \Gamma_j^{\epsilon}\} = \emptyset$, we have $I_{i,j}^{\tbf} = 0$. Otherwise, 
\begin{equation*}
\tau := \inf \{t_2 \in \R_{+}: \Psi_0^{(t_2, \tbf)}(x) \in \Gamma_j^{\epsilon}\}
\end{equation*}
is a number in $[0, \infty)$ such that 
\begin{equation*}
\Psi_0^{(t_2, \tbf)}(x) \notin \Gamma_j^{\epsilon}, \quad t_2 < \tau.
\end{equation*}
The estimate in~\eqref{eq:exit_time_estimate} then implies  
\begin{equation*}
\left \lvert \det U \left (\Psi_0^{(\tau+t, \tbf)}(x) \right) \right \rvert \geq \eps, \quad t \geq \tau_{\epsilon}, 
\end{equation*}
so for every $t \geq \tau_{\eps}$, we have $\Psi_0^{(\tau + t, \tbf)}(x) \notin \Gamma^{\eps}_j$. As a result, the right side of~\eqref{eq:case_2_est} is less than 
\begin{equation*}
\frac{1}{\lambda_{i_{n+1}}} \int_{\tau}^{\tau + \tau_{\epsilon}} e^{- \lambda_{i_n} t_2} \ dt_2 \leq \frac{\tau_{\epsilon}}{\lambda_{i_{n+1}}}.  
\end{equation*}
It follows that 
\begin{equation*}
\M_{n+2}^{\eps, \lambda_0, \lambda_1}(x) \leq \frac{4 \tau_{\epsilon}}{\lambda_0 \wedge \lambda_1} \M_n^{\eps, \lambda_0, \lambda_1}(x).  
\end{equation*} 
\epf 

\bigskip

\subsection{Proof of Theorem~\ref{thm:bounded}, part (1), for $\lambda_1 > \alpha + \beta$}     \label{ssec:proof_bounded_fast_switching}

In this subsection, we prove part (1) of Theorem~\ref{thm:bounded} in the case where $\lambda_1 > \alpha + \beta$.  Under this additional assumption, we can give a simpler proof than in the general case treated in Subsection~\ref{ssec:proof_bounded_interior}. In light of Lemma~\ref{lm:Hennion} and Corollary~\ref{co:convergence_to_0}, it is enough to show that there is $c>0$ such that for small $\epsilon > 0$,  
\begin{equation}   \label{eq:A_series}
\sum_{k=1}^{\infty} \AC_k^{\epsilon} \rho_{i_{k+1}}(x) < c, \quad x \in \Gamma^{\circ}. 
\end{equation}
Let 
$$
m = \max\{(\lambda_0 - (\alpha + \beta))^{-1}, (\lambda_0 - (\alpha + \beta))^{-1} (\lambda_1 - (\alpha + \beta))^{-1}\}.
$$ 
It is easy to see that 
$$ 
\max \left\{\M_1^{\eps, \lambda_0 -(\alpha + \beta), \lambda_1 - (\alpha+\beta)}(x), \M_2^{\eps, \lambda_0 - (\alpha + \beta), \lambda_1 - (\alpha +\beta)}(x) \right\} \leq m
$$ 
for all $x \in \Gamma^{\circ}$. Lemma~\ref{lm:geometric_series} implies 
\begin{equation}    \label{eq:geometric_estimate}
\M_k^{\eps, \lambda_0 - (\alpha+\beta), \lambda_1 - (\alpha +\beta)}(x) \leq f(\epsilon)^{\frac{k - 2}{2}} m = f(\epsilon)^{-1} m f(\epsilon)^{\frac{k}{2}}, \quad k \in \N. 
\end{equation}
Combining~\eqref{eq:geometric_estimate} with Lemmas~\ref{lm:change_of_variables_1} and~\ref{lm:estimates}, we obtain 
\begin{equation*}   
\sum_{k=1}^{\infty} \AC_k^{\eps} \rho_{i_{k+1}}(x) \leq \frac{2 \lambda_0 \lambda_1}{\epsilon} \biggl(1+ f(\epsilon)^{-1} m  \sum_{k=1}^{\infty} \lambda_0(k) f(\epsilon)^{\frac{k}{2}} \biggr), \quad x \in \Gamma^{\circ}. 
\end{equation*}
If $\epsilon$ is so small that $(\lambda_0 \vee \lambda_1) \sqrt{f(\epsilon)} < 1$, the series on the right converges. This completes the proof.
\bigskip

\subsection{Proof of Theorem~\ref{thm:bounded}, part (1)}       \label{ssec:proof_bounded_interior}

Let $\lambda_0 > \alpha + \beta$ and $\lambda_1 \in (\beta, \alpha + \beta]$ (recall that the case $\lambda_1 > \alpha + \beta$ has been taken care of in Subsection~\ref{ssec:proof_bounded_fast_switching}). For $n \in \N$, $\eps > 0$, $x \in \Gamma^{\circ}$, and $\lambda_0, \lambda_1 \in \R$, set 
\begin{equation}     \label{eq:widehat_M} 
\widehat \M_n^{\eps, \lambda_0, \lambda_1}(x) = \int_{M_n^{\eps}(x)} e^{-\langle \lambda_0^{(n+1)}, (\tau_{i_n}(\Psi_0^{\tbf} x), \tbf) \rangle} \ d \tbf,  
\end{equation}
where one should recall that $\tau_i(x)$ was defined in the statement of Lemma~\ref{lm:estimates}. We have 
\begin{equation}     \label{eq:A_series_estimate}
\sum_{k=1}^{\infty} \AC_k^{\eps} \rho_{i_{k+1}}(x) \leq \frac{2 \lambda_0 \lambda_1}{\eps} \biggl(1 + \sum_{k=1}^{\infty} \lambda_0(k) \widehat \M_k^{\eps, \lambda_0 - (\alpha + \beta), \lambda_1 - (\alpha + \beta)}(x) \biggr). 
\end{equation}
Now we need to estimate $\widehat \M_n^{\eps, \lambda_0, \lambda_1}$ for $\lambda_0 > 0$ and $\lambda_1 \in (-\alpha, 0]$.   

\begin{lemma}    \label{lm:geometric_series_general}
There are functions $c(\eps)$ and $f(\eps)$ such that $\lim_{\eps \downarrow 0} f(\eps) = 0$ and 
\begin{equation}    \label{eq:gmg}
\widehat \M_n^{\eps, \lambda_0, \lambda_1}(x) \leq c(\eps) f(\eps)^n, \quad n \in \N, \ x \in \Gamma^{\circ}, \ \lambda_0 > 0, \ \lambda_1 \in (-\alpha, 0]. 
\end{equation}
\end{lemma} 

By Lemma~\ref{lm:geometric_series_general}, the right side of~\eqref{eq:A_series_estimate} is bounded by 
$$
\frac{2 \lambda_0 \lambda_1}{\eps} \biggl(1 + c(\eps) \sum_{k=1}^{\infty} (\lambda_0 f(\eps))^k \biggr),  
$$
which does not depend on $x$ and is finite for $\eps$ so small that $\lambda_0 f(\eps) < 1$. To complete the proof of Theorem~\ref{thm:bounded}, part (1), it remains to show Lemma~\ref{lm:geometric_series_general}. We will do this at the end of this subsection. 

For $i \in \{0,1\}$, $\eps > 0$ and $x \in \Gamma^{\circ}$ such that $\lvert \det U(x) \rvert < \eps$, let 
$$ 
\tau_i^{\eps}(x) = \sup \{t \geq 0: \lvert \det U(\Psi_i^t x) \rvert < \eps\} \wedge \sup \{t \geq 0: \Psi_i^t(x) \in \Gamma^{\circ}\}. 
$$ 
It is easy to see that $\tau_1^{\eps}(x) = \tau_0^{\eps}((1,1)-x)$. 
For integers $j \geq 2$, $i \in \{0,1\}$, and $\eps > 0$, we define 
$$ 
V_i^{\eps}(j) = \{x \in \Gamma^{\circ}: \lvert \det U(x) \rvert < \eps, \ \tau_i^{\eps}(x) \in (j-1, j]\}. 
$$ 
In order to deal with short exit times for the strip of points $x \in \Gamma^{\circ}$ such that $\lvert \det U(x) \rvert < \eps$, we also define 
\begin{align*}
V_i^{\eps}(1) =& \{x \in \Gamma^{\circ}: \lvert \det U(x) \rvert < \eps, \ \tau_i^{\eps}(x) \in (\eps^{\frac{\beta}{\alpha}}, 1]\}, \\
V_i^{\eps}(0) =& \{x \in \Gamma^{\circ}: \lvert \det U(x) \rvert < \eps, \ \tau_i^{\eps}(x) \in (0,\eps^{\frac{\beta}{\alpha}}]\}. 
\end{align*}
In the following lemma, we analyze the interplay of exit times $\tau_0^{\eps}(x)$ and $\tau_1^{\eps}(x)$. 

\begin{lemma}  \label{lm:integral_estimate}
There exists a family of constants $(m^{\eps}(j))_{\eps > 0, j \in \N}$ with the following properties. 
\begin{enumerate} 
\item For $\eps > 0$ sufficiently small and for any $j \in \N$, one has $\tau_i^{\eps}(x) < m^{\eps}(j)$ for every $i \in \{0,1\}$ and $x \in V^{\eps}_{1-i}(j)$. 
\item One has  
$$ 
\lim_{\eps \downarrow 0} \sum_{j=1}^{\infty} j e^{\nu j} m^{\eps}(j) = 0, \quad \nu \in [0, \alpha). 
$$ 
\end{enumerate} 
\end{lemma}

\bpf We first define the constants $(m^{\eps}(j))_{\eps > 0, j \in \N}$ and then verify the asserted properties (1) and (2). In a first step, we show that there is $c > 0$ such that for every $x \in \Gamma^{\circ}$ with $x_1 < \tfrac{1}{2}$, 
\begin{equation}     \label{eq:mvt}
\tau_0^{\eps}((1,1) - x) < c x_1. 
\end{equation}
For $y \in \Gamma^{\circ}$, we have 
$$ 
(1-y_1)^{\frac{\beta}{\alpha}} > 1 - y_2
$$ 
and there is a unique $t^*(y) > 0$, given by 
$$ 
t^*(y) = \sup \{t \geq 0: \Psi_0^t(y) \in \Gamma^{\circ}\}, 
$$ 
such that 
$$ 
\left(1 - e^{\alpha t^*(y)} y_1 \right)^{\frac{\beta}{\alpha}} = 1 - e^{\beta t^*(y)} y_2. 
$$ 
For fixed $y_1 \in (0,1)$, one observes that $t^*(y)$ is an increasing function of $y_2$: the larger $y_2$, the smaller is $1-y_2$, and the larger is the absolute value of the derivative
$$ 
\frac{d}{dt} \left(1 - e^{\beta t} y_2 \right) = -\beta e^{\beta t} y_2, 
$$ 
i.e. the longer it takes the decreasing term $(1-e^{\alpha t} y_1)^{\frac{\beta}{\alpha}}$ to catch up with $1-e^{\beta t} y_2$. Now, consider the function  
$$ 
a(x_1, x_2, t) = 1 - e^{\alpha t} (1-x_1) - \left \lvert 1 - e^{\beta t} (1-x_2) \right \rvert^{\frac{\alpha}{\beta}}, \quad (x_1, x_2, t) \in \R^3. 
$$ 
Since $\tfrac{\alpha}{\beta} > 1$, the function $a$ is $C^1$ on $\R^3$. As $a(0,0,0) = 0$ and $\partial_t a(0,0,0) = -\alpha < 0$, the implicit function theorem implies that there is an open neighborhood $U$ of $(0,0)$ and a $C^1$ function $b : U \to \R$ that is uniquely determined by $b(0,0) = 0$ and 
$$ 
a(x_1, x_2, b(x_1, x_2)) = 0, \quad (x_1, x_2) \in U. 
$$ 
The reason for defining $a$ in the first place is that the function $b$ induced by $a$ satisfies  
$$ 
b(x) = t^*((1,1)-x), \quad x \in \Gamma^{\circ} \cap U. 
$$ 
Thus, for $x \in \Gamma^{\circ}$ sufficiently close to $(0,0)$, we have 
$$ 
\tau_0^{\eps}((1,1)-x) \leq t^*((1,1)-x) \leq b(x_1, 1-(1-x_1)^{\frac{\beta}{\alpha}}). 
$$ 
The rightmost expression in the chain of inequalities above is a $C^1$ function of $x_1$ for $x_1$ close to $0$, and yields $0$ when evaluated at $x_1 = 0$. The mean-value theorem then implies~\eqref{eq:mvt} for a suitable $c > 0$ and for $x \in \Gamma^{\circ}$ with $x_1$ sufficiently close to $0$. As $\tau_0^{\eps}((1,1)-x)$ is bounded for $x \in \Gamma^{\circ}$ such that $x_1 < \tfrac{1}{2}$, we can extend~\eqref{eq:mvt} to the set of such $x$ by choosing a larger $c$.  

Let $q$ be the unique real number greater than $1$ such that $q - q^{\frac{\beta}{\alpha}} = 1$, and let $c' > 0$ be a constant such that 
$$
\frac{1}{e^{(\alpha - \beta) j} - 1} < c' e^{(\beta - \alpha) j}, \quad j \geq 1.
$$ 
Then, define   
$$ 
m^{\eps}(j) = \begin{cases}  
                          c \min \left\{q e^{-\alpha (j-1)}, \eps c' \frac{2}{\alpha \beta} e^{(\beta - \alpha) (j-1)} \right\}, & \quad j \geq 2, \\
c \frac{2}{\alpha \beta (\alpha - \beta)} \eps^{1-\frac{\beta}{\alpha}}, & \quad j=1. 
\end{cases}
$$ 
Fix $\nu \in [0,\alpha)$, and choose $\eta > 0$ such that $\alpha - \beta < \nu + \eta < \alpha$. Let $C > 0$ be so large that $j+1 \leq C e^{\eta j}$ for all $j \in \N$. 
In addition, let $\gamma > 0$ be so small that $\gamma (\nu + \eta + \beta - \alpha) < 1$. As $e^{\nu + \eta - \alpha} < 1$, we obtain for sufficiently small $\eps$ the estimate 
\begin{align*}
& \sum_{j=1}^{\infty} j e^{\nu j} m^{\eps}(j) \\
\leq& c e^{\nu} \biggl(\frac{2}{\alpha \beta (\alpha - \beta)}  \eps^{1-\frac{\beta}{\alpha}} + C \eps c' \frac{2}{\alpha \beta} \sum_{j=1}^{\lfloor -\ln(\eps^{\gamma}) \rfloor} e^{(\eta+\nu -\alpha + \beta) j} + C q \sum_{j = \lfloor -\ln(\eps^{\gamma}) \rfloor +1}^{\infty} e^{(\eta +\nu -\alpha)j} \biggr) \\
\leq&  c e^{\nu} \left(\frac{2}{\alpha \beta (\alpha - \beta)} \eps^{1-\frac{\beta}{\alpha}} - \frac{2 C c' \gamma}{\alpha \beta} \ln(\eps) \eps^{1 - \gamma (\nu + \eta + \beta - \alpha)} + C q \frac{\eps^{\gamma (\alpha - \nu - \eta)}}{1-e^{\nu + \eta  - \alpha}} \right),
\end{align*}
and the right side converges to $0$ as $\eps \downarrow 0$.  This establishes property (2). 

We now show property (1). For symmetry reasons, we can restrict ourselves to the case $i=1$, i.e. we will show that $\tau_1^{\eps}(x) < m^{\eps}(j)$ for all $x \in V_0^{\eps}(j)$.  
Let $x \in V_0^{\eps}(j)$ for some $j \geq 2$. Then we have $\tau_0^{\eps}(x)  > j-1$, which implies 
$$ 
\sup \{t \geq 0: \lvert \det U(\Psi_0^t x) \rvert < \eps \} > j-1. 
$$ 
Therefore, 
$$ 
\alpha \beta \left(e^{\alpha (j-1)} x_1 - e^{\beta (j-1)} x_2 \right) < \eps
$$ 
and, setting $y := e^{\alpha (j-1)} x_1$ and using that $x_1 > x_2^{\frac{\alpha}{\beta}}$,   
\begin{align} 
\frac{\eps}{\alpha \beta} >& e^{\alpha (j-1)} x_1 - e^{\beta (j-1)} x_2 \label{eq:define_m} \\
>& e^{\alpha (j-1)} x_1 - e^{\beta (j-1)} x_1^{\frac{\beta}{\alpha}} = y - y^{\frac{\beta}{\alpha}}. \notag 
\end{align}
As $y \mapsto y - y^{\frac{\beta}{\alpha}}$ is increasing for $y \geq 1$ and negative for $y \in (0,1)$, the definition of $q$ implies $e^{\alpha (j-1)} x_1 < q$ and thus 
$$ 
x_1 < q e^{-\alpha (j-1)}, 
$$ 
provided that $\eps < \alpha \beta$. Since $\lvert \det U(x) \rvert < \eps$ and thus $x_2 < x_1 + \tfrac{\eps}{\alpha \beta}$, the right side of the first line of~\eqref{eq:define_m} is also greater than 
$$ 
e^{\alpha (j-1)} x_1 - e^{\beta (j-1)} x_1 - e^{\beta (j-1)} \frac{\eps}{\alpha \beta},  
$$ 
whence it follows that 
\begin{align*} 
x_1 < \frac{\eps}{\alpha \beta} \frac{1+e^{\beta (j-1)}}{e^{\alpha (j-1)} - e^{\beta (j-1)}} =& \frac{\eps}{\alpha \beta} \frac{1+e^{-\beta (j-1)}}{e^{(\alpha -\beta)(j-1)} - 1} \\
<& \eps \frac{2}{\alpha \beta} \frac{1}{e^{(\alpha - \beta)(j-1)}-1} < \eps c' \frac{2}{\alpha \beta} e^{(\beta - \alpha) (j-1)}. 
\end{align*} 
So far, we have shown that $x_1 < c^{-1} m^{\eps}(j)$. The asserted inequality $\tau_1^{\eps}(x) < m^{\eps}(j)$ then follows from $\tau_1^{\eps}(x) < c x_1$, which is equivalent to~\eqref{eq:mvt}. It remains to consider the case $j=1$. Similarly to the case $j \geq 2$, 
$$ 
x_1 < \frac{2 \eps}{\alpha \beta (e^{(\alpha - \beta) \eps^{\frac{\beta}{\alpha}}} - 1)} \leq c^{-1} m^{\eps}(1). 
$$ 
We can then conclude as in the case $j \geq 2$. 
\epf 

\bigskip

\bpf[Proof of Lemma~\ref{lm:geometric_series_general}] Fix $n \in \N$, $\lambda_0 > 0$, $\lambda_1 \in (-\alpha, 0]$ and $x \in \Gamma^{\circ}$. Since 
$$ 
\{x \in \Gamma^{\circ}: \lvert \det U(x) \rvert < \eps\} = \bigcup_{j=0}^{\infty} V_0^{\eps}(j) = \bigcup_{j=0}^{\infty} V_1^{\eps}(j),
$$ 
and since the sets $(V_i^{\eps}(j))_{j \in \N_0}$ are disjoint for $i \in \{0,1\}$, we have 
\begin{align}   \label{eq:gmg_1}
\widehat \M^{\eps, \lambda_0, \lambda_1}_n(x) =& \sum_{(j_1, \ldots, j_n) \in \N_0^n} \int_{t_n \in \R_{+}: \Psi_0^{t_n}(x) \in V_1^{\eps}(j_n)} \ d t_n \  e^{-\lambda_0 t_n}  \\
& \int_{t_{n-1} \in \R_{+}: \Psi_0^{(t_{n-1}, t_n)}(x) \in V_0^{\eps}(j_{n-1})} \ d t_{n-1} \  e^{-\lambda_1 t_{n-1}} \ldots \notag \\
& \int_{t_1 \in \R_{+}: \Psi_0^{(t_1, \ldots, t_n)}(x) \in V_{i_n}^{\eps}(j_1)}  \ d t_1 e^{- \lambda_{1-i_n} t_1 - \lambda_{i_n} \tau_{i_n}(\Psi_0^{(t_1, \ldots, t_n)} x))},  \notag
\end{align}
where $i_n = 0$ for $n$ even, and $i_n = 1$ for $n$ odd. Fix $(j_1, \ldots, j_n) \in \N_0^n$ and a sequence of switching times $(t_2, \ldots, t_n)$ such that 
$$ 
\Psi_0^{(t_k, \ldots, t_n)}(x) \in V_{i_{n-k+1}}^{\eps}(j_k), \quad 2 \leq k \leq n. 
$$ 
We would like to estimate the integral in the third line of~\eqref{eq:gmg_1}. Set $z := \Psi_0^{(t_2, \ldots, t_n)}(x)$ and let $t_1 \in \R_+$ such that $\Psi_{i_{n+1}}^{t_1}(z) \in V_{i_n}^{\eps}(j_1)$. Since $z \in V_{i_{n-1}}^{\eps}(j_2)$, we have $\lvert \det U(z) \rvert < \eps$. Besides, as $\Psi_{i_{n+1}}^{t_1}(z) \in V_{i_n}^{\eps}(j_1)$, we have $\Psi_{i_{n+1}}^{t_1}(z) \in \Gamma^{\circ}$ and $\lvert \det U(\Psi_{i_{n+1}}^{t_1}(z)) \rvert < \eps$. Hence,   
\begin{equation}      \label{eq:peace}
t_1 < \tau_{i_{n+1}}^{\eps}(z) \leq j_2 + 1. 
\end{equation}
For $j_2 \geq 1$, we even have $\tau_{i_{n+1}}^{\eps}(z) \leq j_2$, but we work with this slightly worse estimate to avoid distinguishing between the cases $j_2 = 0$ and $j_2 \geq 1$.  
Recall that $\tau_i(x)$ was defined in Lemma~\ref{lm:estimates}. We claim that for every $\eps > 0$, 
\begin{equation}   \label{eq:d_claim}  
d(\epsilon) := \sup_{x \in \Gamma^{\circ}: \lvert \det U(x) \rvert < \eps} \left( \tau_0(x) - \tau_0^{\eps}(x) \right) < \infty. 
\end{equation} 
Let $x \in \Gamma^{\circ}$ such that $\lvert \det U(x) \rvert < \eps$, and assume without loss of generality that $\tau_0(x)$ is strictly larger than $\tau_0^{\eps}(x)$. Then, 
$$ 
\tau_0(x) - \tau_0^{\eps}(x) = \tau_0(y(x)), 
$$ 
where 
$$
y(x) = \Psi_0^{\tau_0^{\eps}(x)}(x) \in \{z \in \Gamma^{\circ}:  \det U(z) = \eps\}. 
$$
Since $\tau_0$ is continuous and since $\tau_0(y)$ converges to $0$ as $y$ approaches either one of the endpoints of the line segment $\{z \in \Gamma^{\circ}: \det U(z) = \eps\}$, the claim in~\eqref{eq:d_claim} follows. Together with~\eqref{eq:peace}, the fact that $\Psi_{i_{n+1}}^{t_1}(z) \in V_{i_n}^{\eps}(j_1)$, and the assumption $\lambda_0 > 0 \geq \lambda_1$, one obtains 
$$ 
e^{-\lambda_{1-i_n} t_1 - \lambda_{i_n} \tau_{i_n}(\Psi_0^{(t_1, \ldots, t_n)} x)} \le e^{-\lambda_1 (j_1 + j_2 + d(\eps) + 2)}, 
$$ 
so the integral in the third line of~\eqref{eq:gmg_1} is bounded from above by 
$$
e^{-\lambda_1 (j_1 + j_2 + d(\eps) + 2)} \ \Leb \left(\left \{t_1 \in \R_{+}: \Psi_{i_{n+1}}^{t_1}(z) \in V^{\eps}_{i_n}(j_1) \right\} \right),  
$$
where $\Leb$ denotes the Lebesgue measure. 
Fix a time 
$$
t_1 \in \left\{ t_1 \in \R_+: \Psi_{i_{n+1}}^{t_1}(z) \in V_{i_n}^{\eps}(j_1) \right\}
$$
and define $y := \Psi_{i_{n+1}}^{t_1}(z)$. If $j_1 \geq 1$, Lemma~\ref{lm:integral_estimate} yields $\tau_{i_{n+1}}^{\eps}(y) < m^{\eps}(j_1)$, so 
\begin{equation}    \label{eq:Lebesgue_time_set} 
\Leb \left( \left\{ t_1 \in \R_+: \Psi_{i_{n+1}}^{t_1}(z) \in V_{i_n}^{\eps}(j_1) \right\} \right) \leq m^{\eps}(j_1). 
\end{equation} 
If $j_1 = 0$, the expression on the lefthand side of~\eqref{eq:Lebesgue_time_set} is bounded from above by 
\begin{equation*}
\sup_{y \in V_{i_{n-1}}^{\eps}(j_2)} \tau_{i_{n+1}}^{\eps}(y) = j_2 + \eps^{\frac{\beta}{\alpha}} \id_{j_2 = 0}. 
\end{equation*}
For $i, j \in \N_0$, define 
\begin{equation*}
h(i,j) = e^{-\lambda_1 (j + 1)} \ (m^{\eps}(i) \id_{i \geq 1} + (j + \eps^{\frac{\beta}{\alpha}} \id_{j = 0}) \id_{i=0}). 
\end{equation*}
The integral in the third line of~\eqref{eq:gmg_1} is thus less than 
\begin{equation*}
e^{-\lambda_1 (j_1 + 1 + d(\epsilon))} h(j_1, j_2). 
\end{equation*}
A similar estimate without the factor $e^{-\lambda_1  (j_1 + 1 + d(\epsilon))}$ applies to each of the other integrals in~\eqref{eq:gmg_1}, with the crucial exception of  
\begin{equation}    \label{eq:exception}
\int_{t_n \in \R_{+}:  \Psi_0^{t_n}(x) \in V_1^{\eps}(j_n)} e^{-\lambda_0 t_n} \ d t_n \leq \frac{1}{\lambda_0}.   
\end{equation}
Then, 
\begin{align}     \label{eq:man}  
& \widehat \M_n^{\eps, \lambda_0, \lambda_1}(x) \\
\leq& \frac{e^{-\lambda_1 d(\epsilon)}}{\lambda_0} \sum_{(j_1, \ldots, j_n) \in \N_0^n} e^{-\lambda_1 (j_1 + 1)} \prod_{l=1}^{n-1}  h(j_l,  j_{l+1})  
= \frac{e^{-\lambda_1 d(\epsilon)}}{\lambda_0} \sum_{A \subset \{1, \ldots, n\}} J_A,  \notag
\end{align}
where 
$$
J_A = \sum_{(j_1, \ldots, j_n) \in \N_A} e^{-\lambda_1 (j_1 +1)} \prod_{l=1}^{n-1} h(j_l, j_{l+1}) 
$$
and  
\begin{equation*} 
\N_A = \{(j_1, \ldots, j_n) \in \N_0^n:  j_l > 0 \ \text{iff} \  l \in A\}, \quad A \subset \{1, \ldots, n\}. 
\end{equation*}
Fix a set $A \subset \{1, \ldots, n\}$.  We want to estimate the term $J_A$.  
If $A = \emptyset$, 
\begin{equation}         \label{eq:empty_A} 
J_A = e^{-\lambda_1} h(0,0)^{n-1} = \eps^{-\frac{\beta}{\alpha}} \left(\eps^{\frac{\beta}{\alpha}} e^{-\lambda_1} \right)^n. 
\end{equation}
Now assume $A \neq \emptyset$. We call a subset $B$ of $A$ a \emph{connected component} 
if $k < l < m$ for $k, m \in B$ and $l \in \{1, \ldots, n\}$ implies $l \in B$ and if no subset of $A$ that strictly contains $B$ has this property. The set $A$ can be written as the disjoint union of its \CCs, and the number of {\CCs} of the set $A$ ranges from $1$ to the cardinality of $A$. We call the number of indices between two adjacent {\CCs} $B_1$ and $B_2$ the \emph{gap} between $B_1$ and $B_2$.  

Let $B_1, \ldots, B_m$ be the {\CCs} of a nonempty set $A \subset \{1, \ldots, n\}$, written in increasing order, i.e. $B_1$ contains the smallest index in $A$, $B_2$ is adjacent to $B_1$, etc.  We denote the sizes of $B_1, \ldots, B_m$ by $s_1, \ldots, s_m$, respectively, and write $B_l = \{j^l_1, \ldots, j^l_{s_l}\}$ for $1 \leq l \leq m$.  For $1 \leq l \leq m - 1$, let $g_l$ denote the gap between $B_l$ and $B_{l+1}$. Let $g_0$ be the number of indices preceding $B_1$ and let $g_m$ be the number of indices following $B_m$.   
Then, $\vartheta := (g_0 - 1)^{+} + \sum_{l=1}^{m-1} (g_l - 1) + (g_m - 1)^{+}$ is the number of instances in which two subsequent indices are both not in $A$. Depending on whether $g_0 > 0$ or $g_m > 0$, there are four cases to consider. We only present the case $g_0, g_m > 0$, and thus omit the cases where $g_0 = 0$ or $g_m = 0$. We have  
\begin{align}
J_A = e^{-\lambda_1} h(0,0)^{\vartheta} \prod_{l=1}^m & \biggl(\sum_{j^l_1 = 1}^{\infty}  h(0,j^l_1)  \label{eq:gmg_4} \\
& \sum_{j^l_2 = 1}^{\infty} h(j^l_1, j^l_2) \ldots \sum_{j^l_{s_l} = 1}^{\infty} h(j^l_{s_l - 1}, j^l_{s_l}) h(j^l_{s_l}, 0) \biggr).  \notag
\end{align}
For $j^l_1, \ldots, j^l_{s_l} \geq 1$, 
\begin{equation*}
h(0,j^l_1) h(j^l_1, j^l_2) \ldots h(j^l_{s_l - 1}, j^l_{s_l}) h(j^l_{s_l}, 0) = e^{-\lambda_1}  j^l_1 \prod_{i = 1}^{s_l} e^{-\lambda_1 (j^l_i + 1)} m^{\epsilon}(j^l_i). 
\end{equation*}
Thus, the right side of~\eqref{eq:gmg_4} can be written as 
$$
e^{-\lambda_1} \left(\eps^{\frac{\beta}{\alpha}} \ e^{-\lambda_1} \right)^{\vartheta}  e^{-\lambda_1 m} \prod_{l=1}^m \biggl( e^{-\lambda_1 s_l} \biggl(\sum_{j^l_1 = 1}^{\infty} j^l_1 e^{-\lambda_1 j^l_1} m^{\epsilon}(j^l_1) \biggr) \biggl( \sum_{j = 1}^{\infty} e^{-\lambda_1 j} m^{\epsilon}(j) \biggr)^{s_l - 1} \biggr).
$$ 
It follows that
\begin{equation}    \label{eq:gmg_6}
J_A \leq e^{-\lambda_1} \left(\eps^{\frac{\beta}{\alpha}} \ e^{-\lambda_1} \right)^{\vartheta} e^{-\lambda_1 (m + \lvert A \rvert)} b(\eps)^{\lvert A \rvert}, 
\end{equation}
where 
$$
b(\eps) = \sum_{j=1}^{\infty} j e^{-\lambda_1 j} m^{\eps}(j),  
$$ 
and where $\lvert A \rvert = \sum_{l=1}^m s_l$ is the cardinality of $A$.  Since $\lambda_1 \in (-\alpha, 0]$, Lemma~\ref{lm:integral_estimate} implies that $\lim_{\eps \downarrow 0} b(\eps) = 0$. If we set 
\begin{equation*}
\tilde f(\epsilon) = \left(\eps^{\frac{\beta}{\alpha}}  \ e^{-\lambda_1} \right) \vee b(\eps), 
\end{equation*}
we obtain 
\begin{equation}   \label{eq:gmg_8}
J_A \leq e^{-\lambda_1} \left( e^{-2 \lambda_1} \tilde f(\epsilon)^{\frac{1}{3}} \right)^n 
\end{equation}
because $\vartheta + \lvert A \rvert \geq \tfrac{n-1}{2}$ and $m + \lvert A \rvert \leq 2n$. 
This completes the estimate of $J_A$ in the case $g_0, g_n > 0$. In each of the remaining three cases, one can show without much effort that $J_A$ is less than or equal to the expression on the right side of~\eqref{eq:gmg_6}. 

From~\eqref{eq:man},~\eqref{eq:empty_A} and~\eqref{eq:gmg_8}, we infer that 
$$
\widehat \M_n^{\eps, \lambda_0, \lambda_1}(x) \leq c(\eps) f(\eps)^n 
$$
for 
$$
c(\eps) = \frac{e^{-\lambda_1 d(\epsilon)}}{\lambda_0}  \eps^{-\frac{\beta}{\alpha}} e^{-\lambda_1}, \quad f(\eps) = 2 e^{-2 \lambda_1} \left(\eps^{\frac{\beta}{\alpha}} \vee \tilde f(\epsilon)^{\frac{1}{3}} \right). 
$$
This completes the proof of Lemma~\ref{lm:geometric_series_general} and thus of Theorem~\ref{thm:bounded}, part (1). 
\epf 

\bigskip

\subsection{Proof of Theorem~\ref{thm:bounded}, part (2)}     

The proof of part (2) is quite similar to the proof of part (1) for $\lambda_1 \in (\beta, \alpha + \beta]$. Let $K$ be a compact subset of $\Gamma_0$, and recall that $\widehat M$ was defined in~\eqref{eq:widehat_M}. Following the proof of part (1), all we need to show is the following lemma. 

\begin{lemma}      \label{lm:geometric_series_general_alt}
There are functions $c(\eps)$ and $f(\eps)$ such that $\lim_{\eps \downarrow 0} f(\eps) = 0$ and 
\begin{equation}    \label{eq:gmg_alt}
\widehat \M_n^{\eps, \lambda_0, \lambda_1}(x) \leq c(\eps) f(\eps)^n, \quad n \in \N, \ x \in K \cap \Gamma^{\circ}, \ \lambda_0, \lambda_1 \in (-\alpha, 0]. 
\end{equation}
\end{lemma}

\bpf Fix $n \in \N$, $\lambda_0, \lambda_1 \in (-\alpha, 0]$ and $x \in K \cap \Gamma^{\circ}$. As in the proof of Lemma~\ref{lm:geometric_series_general}, we use the representation for $\widehat \M_n^{\eps, \lambda_0, \lambda_1}(x)$ in~\eqref{eq:gmg_1}. If we let $\lambda = \lambda_0 \wedge \lambda_1$, the integral in the third line of~\eqref{eq:gmg_1} is less than 
$$ 
e^{-\lambda (j_1 + 1 + d(\eps))} h(j_1, j_2), 
$$ 
where 
$$ 
h(i,j) = e^{-\lambda (j+1)} \left(m^{\eps}(i) \id_{i \geq 1} + \left(j + \eps^{\frac{\beta}{\alpha}} \id_{j=0} \right) \id_{i=0} \right). 
$$
As in the proof of Lemma~\ref{lm:geometric_series_general}, there are similar estimates for the other integrals in~\eqref{eq:gmg_1}, but now the integral on the left side of~\eqref{eq:exception} is less than 
$$ 
h_K(j_n) := e^{-\lambda (T_K +1)} \left(m^{\eps}(j_n) \id_{j_n \geq 1} + T_K \id_{j_n = 0} \right). 
$$ 
Here, $T_K$ is a positive integer that depends only on the compact set $K$, and whose existence follows from the assumption $(0,0) \notin K$.  The proof can then essentially be completed as the one of Lemma~\ref{lm:geometric_series_general}, with $\tfrac{1}{\lambda_0}$ replaced by $h_K(j_n)$. 
\epf 

\bigskip

\subsection{Proof of Theorem~\ref{thm:bounded}, part (3)}    \label{sec:proof_bounded_a_l}

The proof strategy is similar to the one from the proof of part (1) for $\lambda_1 > \alpha + \beta$.  Here, however, we need to make sure that the compositions of backward trajectories do not become arbitrarily close to the critical points of $u_0$ and $u_1$. This is accomplished by letting the width of the strip around the diagonal shrink to zero as we move backward in time. The procedure only works because we require $K$ to be a positive distance away from the boundary curve $\partial \Gamma_0$.   
For $n \in \N$, $x \in \Gamma^{\circ}$, and $\eps > 0$, let 
\begin{align*} 
\sigma M_n^{\eps}(x) =& \left\{\tbf \in T_0^n(x): \left \lvert \det U(\Psi_0^{(t_{n-j}, \ldots, t_n)} x) \right \rvert < \eps 2^{-j}, \ 0 \leq j \leq n-1 \right\}, \\
\sigma S_n^{\eps}(x) =& \left\{(t_1, t_2, \tbf) \in T_0^{n+1}(x): \tbf  \in \sigma M_{n-1}^{\eps}(x), \ \left \lvert \det U(\Psi_0^{(t_2, \tbf)} x) \right \rvert > \eps 2^{-(n-1)} \right\}.  
\end{align*}
For $h \in L^1(\Gamma^{\circ})$, define 
\begin{align*}
\sigma \AC_n^{\eps} h(x) =& \int_{\sigma S_n^{\eps}(x)} \lambda_0(n+1) e^{- \langle \lambda_0^{(n+1)}, \tbf \rangle} J_0^{\tbf}(x) h(\Psi_0^{\tbf} x) \ d \tbf, \\
\sigma \Bc_n^{\eps} h(x) =& \int_{\sigma M_n^{\eps}(x)} \lambda_0(n) e^{- \langle \lambda_0^{(n)}, \tbf \rangle} J_0^{\tbf}(x) h(\Psi_0^{\tbf} x) \ d \tbf. 
\end{align*}
The letter $\sigma$ stands for ``shrinking'' (referring to the strip around the diagonal) and is meant to help distinguish the notation from the one introduced at the beginning of Subsection~\ref{ssec:diagonal}. 

We need analogs of Lemmas~\ref{lm:Hennion} and~\ref{lm:change_of_variables_1}, which we state in Lemma~\ref{lm:Hennion_2} below. The proofs of these modified statements are almost identical to the proofs of the original ones, and we omit them.  

\begin{lemma}     \label{lm:Hennion_2}
For any $x \in \Gamma^{\circ}$ and $\epsilon > 0$, the following statements hold.
\begin{enumerate}
\item For any $n \in \N$,
$$
\rho_0(x) = \sum_{k=1}^n \sigma \AC_k^{\eps} \rho_{i_{k+1}}(x) + \sigma \Bc_n^{\eps} \rho_{i_n}(x) = \sum_{k=1}^{\infty} \sigma \AC_k^{\eps} \rho_{i_{k+1}}(x).  
$$
\item  We have   
$$
\sigma \AC_1^{\eps} \rho_{i_2}(x) = \Ic_0^{\eps}(x), 
$$
and for any $k > 1$,  
\begin{equation*} 
\sigma \AC_k^{\eps} \rho_{i_{k+1}}(x) = \int_{\sigma M_{k-1}^{\eps}(x)} \lambda_0(k-1) e^{ \langle (\alpha + \beta) \mathbbm{1} - \lambda_0^{(k-1)}, \tbf \rangle} \Ic_{i_{k+1}}^{\eps 2^{-(k-1)}}(\Psi_0^{\tbf} x) \ d \tbf.
\end{equation*}
\end{enumerate} 
\end{lemma} 

Next, we formulate an analog of Lemma~\ref{lm:geometric_series}. For $n \in \N$, $\epsilon > 0$, $x \in \Gamma^{\circ}$, and $\lambda_0, \lambda_1 \in \R$, set 
\begin{equation*}
\sigma \M_n^{\eps, \lambda_0, \lambda_1}(x) = \int_{\sigma M_n^{\eps}(x)} e^{- \langle \lambda_0^{(n)}, \tbf \rangle} \ d \tbf. 
\end{equation*} 

\begin{lemma}    \label{lm:geometric_series_2}
For any $\lambda_0, \lambda_1 \in \R$, there is a function $f(\epsilon)$ such that $\lim_{\epsilon \downarrow 0} f(\epsilon) = 0$ and 
\begin{equation*}
\sigma \M_{n+1}^{\eps, \lambda_0, \lambda_1}(x) \leq f(\epsilon) \sigma \M_n^{\eps, \lambda_0, \lambda_1}(x), \quad n \in \N, \ x \in K \cap \Gamma^{\circ}. 
\end{equation*}
\end{lemma}

Before proving Lemma~\ref{lm:geometric_series_2}, we carry out some preliminary work. 
For $\delta \in (0, \tfrac{1}{2})$, let 
\begin{equation*}     
\Gamma(\delta) = \{x \in \Gamma^{\circ}: \delta < x_2 < 1 - \delta\}. 
\end{equation*}

\begin{lemma}    \label{lm:delta} 
There is $\delta > 0$ such that for $\epsilon > 0$ sufficiently small, 
\begin{equation*}
\Psi_0^{\tbf}(x) \in \Gamma(\delta), \quad x \in K \cap \Gamma^{\circ}, \ n \in \N, \ \tbf \in \sigma M_n^{\eps}(x). 
\end{equation*}
\end{lemma} 

\bpf We will at times use the notation $[x]_2$ for the second component of a point $x \in \R^2$. For $\delta > 0$ and $n \in \N$, we set 
$$
\mathcal{S}(\delta, n) := 2 \delta - \frac{\delta}{2} \sum_{k=1}^{n} 2^{-k+1} 
$$
to simplify notation. 
Due to symmetries and an induction argument, it is enough to show that there is $\delta > 0$ such that for $\eps > 0$ sufficiently small, the following statements hold: 
\begin{enumerate} 
\item For any $x \in K \cap \Gamma^{\circ}$ and $t \in \sigma M_1^{\eps}(x)$, one has 
$$
\Psi_0^t(x) \in \Gamma(2 \delta); 
$$
\item For any even positive integer $n$, $y \in \Gamma(\mathcal{S}(\delta, n))$ such that $\lvert \det U(y) \rvert < \eps 2^{-(n-1)}$, and $t > 0$ such that $\Psi_0^t(y) \in \Gamma^{\circ}$ and $\lvert \det U(\Psi_0^t y) \rvert < \eps 2^{-n}$, one has 
$$
\Psi_0^t(y) \in \Gamma(\mathcal{S}(\delta, n+1)).  
$$
\end{enumerate}
First we specify $\delta$.  For $x \in \R^2$ and $t \in \R$, consider the function 
$$
\delta(x) = x_2 \vee (1-x_2). 
$$
According to Lemma~\ref{lm:basics}, for any $i \in \{0,1\}$ and $x \in \Gamma^{\circ}$, there is a unique $\theta_i(x) \in \R$ such that $\det U(\Psi_i^{\theta_i(x)} x) = 0$ and $\Psi_i^{\theta_i(x)}(x) \in \Gamma^{\circ}$.  We define 
$$
\zeta = \sup_{x \in K \cap \Gamma^{\circ}} \delta \left( \Psi_0^{\theta_0(x)} x \right), 
$$ 
which is strictly less than $1$ because $K$ is compact and does not intersect $\partial \Gamma_0$.  Then, we set 
\begin{equation*}
\delta := \frac{1-\zeta}{3} \wedge \frac{1}{2} \inf_{x \in K} x_2.  
\end{equation*}
Let $\tilde \eps > 0$ be so small that the closure of 
$$ 
\Xi(\delta, \tilde \eps) := \{x \in \Gamma_{\delta}: \lvert \det U(x) \rvert < \tilde \eps\}
$$ 
is contained in $\Gamma^{\circ}$.  Set
$$ 
\vartheta = \inf_{x \in (K \cap \Gamma^{\circ}) \cup \Xi(\delta, \tilde \eps)} [\Psi_0^{\theta_0(x)}(x)]_1 > 0.
$$
For any $x \in \R^2$ such that $\delta \leq x_1 = x_2 \leq 1$, the formula in part (2) of Lemma~\ref{lm:basics} implies 
$$ 
\frac{d}{ds} \det U(\Psi_0^s x)  \vert_{s=0} \geq \alpha \beta (\alpha - \beta) \delta,
$$
so by a compactness argument there is $r > 0$ such that 
\begin{equation}      \label{eq:petals} 
\frac{d}{ds} \det U(\Psi_0^s x)\vert_{s=0} \geq \frac{\alpha \beta}{2} (\alpha - \beta) \delta
\end{equation}   
for every $x \in \R^2$ such that $\delta \leq x_2 \leq 1$ and $\lvert \det U(x) \rvert < r$. 
Next, observe that since $u_0$ and $u_1$ are bounded on the compact set $\Gamma$, there is $C > 0$ such that 
\begin{equation*}
\lvert \partial_t \Psi_i^t(x) \rvert = \lvert u_i(\Psi_i^t(x)) \rvert \leq C 
\end{equation*}
for every $i \in \{0,1\}$, $x \in \Gamma^{\circ}$, and $t \geq 0$ for which $\Psi_i^t(x) \in \Gamma^{\circ}$. 

We proceed to the proof of statements (1) and (2). We will assume that $\eps$ is sufficiently small with respect to $\delta$, $\tilde \eps$, $\vartheta$, and $C$ for the estimates given above to hold.    
First we prove statement (1), which plays the role of the base case in an induction argument. Fix a point $x \in K \cap \Gamma^{\circ}$. For any $t \in \sigma M_1^{\eps}(x)$, one has  
$$ 
2 \delta \leq x_2 < [\Psi_0^t(x)]_2, 
$$ 
because $t \mapsto [\Psi_0^t(x)]_2$ is increasing. To show the estimate 
$$ 
1 - 2 \delta > [\Psi_0^t(x)]_2, 
$$ 
assume without loss of generality that $t > \theta_0(x)$.  Set $c(s) = (c_1(s), c_2(s)) := \Psi_0^s(x)$ for $s \in [\theta_0(x),t]$.  Then, 
\begin{equation}     \label{eq:lower_eps}
\epsilon > \det U(c(t)) = \det U(c(t)) - \det U(c(\theta_0(x))) = \int_{\theta_0(x)}^t \frac{d}{ds} \det U(c(s)) \ ds.
\end{equation} 
For $s \in [\theta_0(x),t]$, one has, again by part (2) of Lemma~\ref{lm:basics}, 
$$
\frac{d}{ds} \det U(c(s)) \geq \alpha \beta (\alpha - \beta) c_1(\theta_0(x)) \geq \alpha \beta (\alpha - \beta) \vartheta.
$$
Together with~\eqref{eq:lower_eps}, this yields 
\begin{equation*}
t - \theta_0(x) < \frac{\epsilon}{\alpha \beta (\alpha - \beta) \vartheta}. 
\end{equation*}
By the mean-value theorem, there is $s^* \in (\theta_0(x),t)$ such that 
\begin{equation*}
c_2(t) - c_2(\theta_0(x)) = (t-\theta_0(x)) c_2'(s^*). 
\end{equation*}
Then, 
\begin{align*}
c_2(t) =& c_2(\theta_0(x)) + c_2(t) - c_2(\theta_0(x)) \\
\leq& \zeta + (t-\theta_0(x)) c_2'(s^*) \leq \zeta + \frac{\epsilon C}{\alpha \beta (\alpha- \beta) \vartheta} \leq 1-3 \delta + \frac{\eps C}{\alpha \beta (\alpha - \beta) \vartheta} < 1 - 2 \delta,
\end{align*}
which completes the proof of statement (1). 

\bigskip

We proceed to the proof of statement (2).  Let $n$ be an even positive integer, let $y \in \Gamma(\mathcal{S}(\delta, n))$
such that $\lvert \det U(y) \rvert < \epsilon 2^{-(n-1)}$, and let $t > 0$ such that $\Psi_0^t(y) \in \Gamma^{\circ}$ and $\lvert \det U(\Psi_0^t y) \rvert < \epsilon 2^{-n}$.  Then, 
\begin{equation*}
[\Psi_0^t(y)]_2 > y_2 > \mathcal{S}(\delta, n)  > \mathcal{S}(\delta, n+1). 
\end{equation*}
It remains to show 
\begin{equation}     \label{eq:psi_lower_bd}
[\Psi_0^t(y)]_2 < 1 - \mathcal{S}(\delta, n+1). 
\end{equation}
As in the proof of statement (1), there is no loss of generality in assuming $t > \theta_0(y)$. Set $d(s) = (d_1(s), d_2(s)) := \Psi_0^s(y)$ for $s \in [\theta_0(y),t]$. Since $y \in \Xi(\delta, \tilde \eps)$, one has $d_1(\theta_0(y)) \geq \vartheta$. Thus, we can essentially repeat the argument from the proof of statement (1) to obtain
\begin{equation*}
t - \theta_0(y) < \frac{\epsilon 2^{-n}}{\alpha \beta (\alpha - \beta) \vartheta}
\end{equation*}  
and 
$$ 
d_2(t) \leq d_2(\theta_0(y)) + \frac{\eps 2^{-n} C}{\alpha \beta (\alpha - \beta) \vartheta}. 
$$
The next step consists in estimating $d_2(\theta_0(y))$ from above. Assume without loss of generality that $d_2(\theta_0(y)) > y_2$. Then,  
\begin{equation}    \label{eq:before_diag}
\epsilon 2^{-(n-1)} > \int_0^{\theta_0(y)} \frac{d}{ds} \det U(d(s)) \ ds 
\end{equation}
in analogy to~\eqref{eq:lower_eps}.  Next, we estimate $\tfrac{d}{ds} \det U(d(s))$ from below. For fixed $s \in [0,\theta_0(y)]$, we claim that 
\begin{equation*}
\det U(d(s)) > -r, 
\end{equation*}
where one should recall that $r$ was introduced in relation to~\eqref{eq:petals}. 
Suppose the claim doesn't hold. Since $\det U(d(0)) > -r$, there is $s^* \in (0,s]$ such that $\det U(d(s^*)) = -r$ and $\det U(d(t)) > -r$ for every $t \in [0,s^*)$. 
Then, 
\begin{align*}
-r =& \det U(d(s^*)) = \det U(y) + \det U(d(s^*)) - \det U(d(0)) \\
>& -\epsilon 2^{-(n-1)} + \int_0^{s^*} \frac{d}{ds} \det U(d(s)) \ ds \geq -\epsilon + s^* \frac{\alpha \beta}{2} (\alpha - \beta) \delta > -r,
\end{align*}
a contradiction.
As a result, the integral on the righthand side of~\eqref{eq:before_diag} is bounded from below by  
$$
\theta_0(y) \frac{\alpha \beta}{2} (\alpha - \beta) \delta.
$$  
Hence, 
\begin{equation*}
\theta_0(y) < \frac{\epsilon 2^{-(n-2)}}{\alpha \beta (\alpha - \beta) \delta},  
\end{equation*}
and we obtain the estimate
\begin{equation*}
d_2(\theta_0(y)) = y_2 + d_2(\theta_0(y)) - d_2(0) < 1 - \mathcal{S}(\delta, n) + \frac{ \epsilon 2^{-(n-2)}}{\alpha \beta (\alpha - \beta) \delta} C.   
\end{equation*}
This yields 
\begin{equation*}
d_2(t) < 1 - \mathcal{S}(\delta, n)  + \left(\frac{4}{\delta} + \frac{1}{\vartheta} \right) \frac{\eps 2^{-n} C}{\alpha \beta (\alpha - \beta)} \leq 1 - \mathcal{S}(\delta, n+1).  
\end{equation*}
\epf 

\bigskip

\bpf [Proof of Lemma~\ref{lm:geometric_series_2}] Let $n \in \N$, $x \in K \cap \Gamma^{\circ}$, $\eps > 0$, and $\lambda_0, \lambda_1 \in \R$. Let $\lambda = \lvert \lambda_0 \rvert \vee \lvert \lambda_1 \rvert$.  
As an immediate consequence of Lemma~\ref{lm:delta}, there is $\delta > 0$, independent of $n$ and $x$, such that for $\epsilon > 0$ sufficiently small, 
\begin{equation}      \label{eq:M_inclusion}
\Psi_0^{(t_{n+1-j}, \ldots, t_{n+1})}(x) \in \Gamma(\delta), \quad 0 \leq j \leq n, \ \tbf \in \sigma M_{n+1}^{\eps}(x). 
\end{equation} 
Let $\tbf = (t_1, t_2, \ldots, t_{n+1}) \in \sigma M_{n+1}^{\eps}(x)$ and set $y := \Psi_0^{(t_2, \ldots, t_{n+1})}(x)$. As we saw in the proof of Lemma~\ref{lm:delta},   
$$
t_1 = t_1 - \theta_{i_n}(y) + \theta_{i_n}(y) < \frac{\eps 2^{-n}}{\alpha \beta (\alpha - \beta) \vartheta} +  \frac{\eps 2^{-(n-2)}}{\alpha \beta (\alpha - \beta) \delta} \leq c \eps,  
$$
where $c > 0$ is a constant that does not depend on $n$. 
Hence, 
$$
\sigma \M_{n+1}^{\eps, \lambda_0, \lambda_1}(x) \leq \int_{\sigma M_n^{\eps}(x)} d \tbf \ e^{- \langle \lambda_0^{(n)}, \tbf \rangle} \int_0^{c \eps} d t_1 \ e^{\lambda t_1} \leq c \eps e^{\lambda c \eps} \sigma \M_n^{\eps, \lambda_0, \lambda_1}(x).  
$$
As $\lim_{\eps \downarrow 0} c \eps e^{\lambda c \eps} = 0$, this completes the proof. 
\epf 

\bigskip

\bpf [Proof of Theorem~\ref{thm:bounded}, part (3)]  By Lemma~\ref{lm:Hennion_2}, we need to show that 
$$ 
\sup_{x \in K \cap \Gamma^{\circ}} \sum_{k=1}^{\infty} \sigma \AC_k^{\eps} \rho_{i_{k+1}}(x) < \infty. 
$$ 
Again by Lemma~\ref{lm:Hennion_2}, we have for $x \in K \cap \Gamma^{\circ}$  
$$
\sum_{k=1}^{\infty} \sigma \AC_k^{\eps} \rho_{i_{k+1}}(x) = \Ic_0^{\eps}(x) + \sum_{k=2}^{\infty} \int_{\sigma M_{k-1}^{\eps}(x)} \lambda_0(k-1) e^{ \langle (\alpha + \beta) \mathbbm{1} - \lambda_0^{(k-1)}, \tbf \rangle} \Ic_{i_{k+1}}^{\eps 2^{-(k-1)}}(\Psi_0^{\tbf} x) \ d \tbf. 
$$
By Lemma~\ref{lm:estimates}, the righthand side is less than 
\begin{align}     \label{eq:bal_2} 
& \frac{2 \lambda_0 \lambda_1}{\epsilon} e^{(\alpha + \beta - \lambda_0) \tau_0(x)}  \\
& + \sum_{k=2}^{\infty} \frac{2 \lambda_0 \lambda_1}{\epsilon} 2^{k-1} \lambda_0(k-1) \int_{\sigma M_{k-1}^{\eps}(x)} e^{\langle (\alpha + \beta) \mathbbm{1} - \lambda_0^{(k)}, (\tau_{i_{k+1}}(\Psi_0^{\tbf} x), \tbf) \rangle} \ d \tbf \notag.   
\end{align}
For any $k \geq 2$ and $\tbf \in \sigma M_{k-1}^{\eps}(x)$, 
we have $\lvert \det U(\Psi_0^{\tbf} x) \rvert < \eps$ 
and, on account of~\eqref{eq:M_inclusion}, we also have $\Psi_0^{\tbf}(x) \in \Gamma(\delta)$.  Therefore, 
\begin{equation*}
\tau_{i_{k+1}}(\Psi_0^{\tbf} x) \leq c,  
\end{equation*}
where $c$ is a finite constant that does not depend on $\epsilon$. The expression in the second line of~\eqref{eq:bal_2} is thus bounded from above by    
\begin{equation}   \label{eq:bal_3}
e^{(\alpha + \beta) c} \frac{2 \lambda_0 \lambda_1}{\epsilon} \sum_{k=2}^{\infty} 2^{k-1} \lambda_0(k-1) \sigma \M_{k-1}^{\eps, \lambda_0 - (\alpha + \beta), \lambda_1 - (\alpha +\beta)}(x). 
\end{equation}
By Lemma~\ref{lm:geometric_series_2}, for $k \geq 2$, 
\begin{equation*}
\sigma \M_{k-1}^{\eps, \lambda_0 - (\alpha + \beta), \lambda_1 - (\alpha + \beta)}(x) \leq f(\epsilon)^{k-2} \sigma \M_1^{\eps, \lambda_0 - (\alpha + \beta), \lambda_1 - (\alpha + \beta)}(x) \leq f(\epsilon)^{k-2} \hat c e^{(\alpha + \beta) \hat c},  
\end{equation*}
where $\lim_{\epsilon \downarrow 0} f(\epsilon) = 0$ and 
$$ 
\hat c = \sup_{x \in K \cap \Gamma^{\circ}} \sup \{t \geq 0: \det U(\Psi_0^t x) \leq \eps\}.
$$ 
Hence, the expression in~\eqref{eq:bal_3} is bounded from above by  
\begin{equation*}
e^{(\alpha + \beta) (c + \hat c)} \frac{2 \lambda_0 \lambda_1 \hat c}{\epsilon f(\epsilon)} \sum_{k=1}^{\infty} (2 f(\epsilon))^k \lambda_0(k),
\end{equation*}
which doesn't depend on $x$ and is finite for $\eps$ sufficiently small.  
\epf 

\bigskip   
\noindent {\bf Acknowledgments:} JCM thanks the NSF grant DMS-1613337 for
partial support of this work. JCM and SDL thank Michael Reed for many
useful discussions earlier in their work in randomly switched ODEs and
this example in particular. TH gratefully acknowledges support through SNF grant $200021-175728/1$. He thanks Florent Malrieu for helpful suggestions related to this example, and Jean-Baptiste Bardet, Michel Bena\"im, Florent Malrieu, Edouard Strickler, and Pierre-Andr\'e Zitt for useful discussions about invariant densities for this and other two-dimensional examples.  YB is grateful to NSF for partial support of his work via grants DMS-1460595 and DMS-1811444. SDL thanks the NSF grants DMS-1944574 and DMS-1814832 for partial support of this work.

\bigskip

\bibliography{regularity_2}
\bibliographystyle{plain}

\end{document}